\newcommand{\qed}{\hfill \ensuremath{\Box}}
\documentclass[12pt,a4paper]{article}
\usepackage{amssymb}
\usepackage[english,francais]{babel}
\usepackage{graphicx}
\usepackage{amsmath}
\usepackage[T1]{fontenc}
\newtheorem{theorem}{Theorem}[section]
\newtheorem{thm}[theorem]{Th\'{e}or\`{e}me}
\newtheorem{defn}[theorem]{D\'{e}finition}
\newtheorem{lem}[theorem]{Lemme}
\newtheorem{cor}[theorem]{Corollaire}
\newtheorem{prop}[theorem]{Proposition}

\newtheorem{rem}[theorem]{Remarque}

\begin{document}
\title{Les extensions intrins\`{e}ques les plus simples des ensembles $\mathbb{R}$ et $\mathbb{N}$ pour une nouvelle Analyse Non Standard}
\author{\footnote{{\bf thierry.bautier@bretagne.iufm.fr} }Thierry Bautier
\\
{\small \ I.U.F.M. de Bretagne, E.S.P.E., site de Vannes (Ecole interne de l'U.B.O.)}
\\
{\small \ 32 Avenue Roosevelt, 56000 Vannes, France}
}
\date{}
\maketitle

\abstract{Le principal r\'{e}sultat de cette recherche est d'avoir prolong\'{e} de la mani\`{e}re la plus simple les ensembles standard $\mathbb{N}$ et $\mathbb{R}$.
\par Les diff\'{e}rentes extensions de $\mathbb{R}$ sont le plus petit sur-anneau int\`{e}gre $\mathbb{R}_o$, le plus petit sur-ensemble totalement ordonn\'{e} et continu $\overline{\mathbb{R}_o}$ et le plus petit sur-corps complet $\Omega$ de $\mathbb{R}$.
\par L'extension de $\mathbb{N}$ est form\'{e}e d'entiers finis ou infiniment grands mais tous d\'{e}finis \`{a} l'unit\'{e} pr\`{e}s, c'est un mod\`{e}le non standard de l'Arithm\'{e}tique de Peano.
\par Toutes ces extensions de $\mathbb{R}$ et $\mathbb{N}$ sont intrins\`{e}ques, ou n\'{e}cessaires. Elles "ne peuvent pas \^{e}tre autres" \`{a} la diff\'{e}rence d'autres structures num\'{e}riques comme les alg\`{e}bres de Weil $\mathbb{R}[X]/(X^{I})$ qui sont contingentes car elles introduisent le choix arbitraire d'un param\`{e}tre $I$.
\\ \\
On d\'{e}veloppe ici pour la premi\`{e}re fois une nouvelle Analyse Non Standard en prolongeant \`{a} $\mathbb{R}_o$ les fonctions $C^{\infty}$ standard puis on \'{e}tend leurs propri\'{e}t\'{e}s diff\'{e}rentielles et int\'{e}grales aux fonctions r\'{e}guli\`{e}res de $\mathbb{R}_o$ dans $\mathbb{R}_o$ qui sont des s\'{e}ries enti\`{e}res convergentes dans un intervalle d'amplitude standard mais non nulle.
\\
Tous ces r\'{e}sultats ont un point de d\'{e}part historique qui est \`{a} chaque fois pr\'{e}cis\'{e}.

\newpage
\begin{center}
\textbf{Abstract}
\end{center}
\textbf{Infinitely small numbers and infinitely large integers defined to a unit for a new Non Standard Calculus }
\\
The main results of this paper are the constructions, both rigorous and intuitive, of the $\aleph$ intrinsic extension of the set of non negative integers $\mathbb{N}$ and the $\Omega$ (resp. $\mathbb{R}_o$, $\overline{\mathbb{R}_o}$) smallest strict over-field (resp. ring, set) of $\mathbb{R}$ set which is totally ordered and complete (resp. entire, continue).
\par
The first differential and integral elements of a new Non Standard Analysis are given. A new and rigorous proof of the Fundamental Theorem of Analysis is given according to G.W.Leibniz historical intuition's.}

\begin{center}
\textbf{Plan}
\end{center}
Introduction (p. 4)
\\ \\
\textbf{1. De nouveaux nombres "r\'{e}els". Premi\`{e}res propri\'{e}t\'{e}s de l'ensemble $\mathbb{R}_o$}\\ \\
1.1. Propri\'{e}t\'{e}s alg\'{e}briques et ordinales de $\mathbb{R}_o=\mathbb{R}[[X]]$ (7)\\
1.2. Prolongement analytique d'une fonction de classe $C^\infty$ de $\mathbb{R}$ dans $\mathbb{R}$ (8)\\
1.3. "Continuit\'{e}", NS*-continuit\'{e} et "diff\'{e}rentiabilit\'{e}" multiple d'une fonction de $\mathbb{R}_o$ dans $\mathbb{R}_o$ (9)\\
1.4. "Diff\'{e}rentiabilit\'{e}" multiple du prolongement analytique d'une fonction $C^{\infty}$ standard (9)\\
1.5. Deux types de diff\'{e}rentielles pour les prolongements analytiques des fonctions analytiques standard et leurs relations r\'{e}ciproques (11)\\
1.6. Sur le contexte historique de cette recherche (13)\\
1.7. Le co\^{u}t d'une r\'{e}volution non standard (15).
\\ \\
\textbf{2. De nouveaux nombres "entiers" d\'{e}finis \`{a} l'unit\'{e} pr\`{e}s. Propri\'{e}t\'{e}s des ensembles $\mathbb{N}[\Sigma]$ et $\aleph$}\\ \\
2.1. Conditions g\'{e}n\'{e}rales pour construire un ensemble de nombres entiers (18)\\
2.2. Construction formelle de $\mathbb{N}[\Sigma]$ et $\mathbb{R}_o^{1+}$ (18)\\
2.3. Deux mod\`{e}les non standard de l'Arithm\'{e}tique de Peano (20)\\
2.4. Le point de vue de I.Newton sur les nombres entiers infiniment grands (21).
\\ \\
\textbf{3. Une d\'{e}monstration du Th\'{e}or\`{e}me Fondamental de l'Analyse Non Standard}\\ \\
3.1. R\'{e}solution de l'\'{e}quation du premier ordre pour le prolongement analytique d'une fonction analytique standard (23)\\
3.2. G\'{e}n\'{e}ralisation aux fonctions r\'{e}guli\`{e}res ou newtoniennes de $\mathbb{R}_o$ dans $\mathbb{R}_0$ (27)\\
3.3. R\'{e}solution d'une \'{e}quation diff\'{e}rentielle d'ordre quelconque pour des fonctions r\'{e}guli\`{e}res (30).
\\ \\
\textbf{4. Le paradis newtonnien et la r\'{e}habilitation posthume de G.W.Leibniz}\\ \\
4.1. Introduction \`{a} cette \'{e}tude d'Histonique (32)\\
4.2. Deux caract\'{e}ristiques pr\'{e}alable de l'\'{e}difice newtonnien (32)\\
4.3. Une modernisation du calcul des fluxions de Newton (32)\\
4.4. Une comparaison des oeuvres de Newton et Leibniz (34)\\
4.5. La Th\'{e}orie de l'Int\'{e}gration de I.Newton (35)\\
4.6. Aux limites de l'\'{e}difice newtonien (36)\\
4.7. Les faiblesses de l'\'{e}difice leibnizien (37)\\
4.8. La Th\'{e}orie de l'Int\'{e}gration de G.W.Leibniz (37)\\
4.9. Conclusion de cette \'{e}tude d'Histonique (38).
\\ \\
\textbf{5. Principaux r\'{e}sultats sur les ensembles de nombres}\\ \\
5.1. Les extensions intrins\`{e}ques les plus simples de $\mathbb{N}$ et $\mathbb{R}$ (40)\\
5.2. Premi\`{e}res propri\'{e}t\'{e}s du corps totalement ordonn\'{e} $(\Omega,+,\times,\leq)$ (41)\\
5.3. Propri\'{e}t\'{e}s de l'alg\`{e}bre totalement ordonn\'{e}e $(\Omega,+,\cdot,\times,\leq)$ (42)\\
5.4. Propri\'{e}t\'{e}s ordinales de $\Omega$ et de $\overline{\Omega}$ (43)\\
5.5. Une nouvelle caract\'{e}risation des structures $(\Omega,+,\times,\leq)$ et  $(\overline{\Omega},\leq)$ (47).\\ \\
Conclusion (p. 50).\\ \\
R\'{e}f\'{e}rences.

\newpage
\section*{Introduction}
Les Analyses Non Standard de A.Robinson [1,2,3] et de J.H.Conway [4] sont aujourd'hui assez bien connues et la nouvelle Analyse Non Standard qui est ici pr\'{e}sent\'{e}e vise les m\^{e}mes objectifs que ces deux th\'{e}ories math\'{e}matiques :
\begin{quotation}
Compl\'{e}ter la droite num\'{e}rique standard par des \'{e}l\'{e}ments infiniment grands et infiniment petits,  en conservant l'essentiel de ses propri\'{e}t\'{e}s.
\end{quotation}
Le niveau technique de cette troisi\`{e}me th\'{e}orie math\'{e}matique est \emph{tr\`{e}s} simple, beaucoup plus simple que celui des deux autres th\'{e}ories et l'on doit s'en expliquer dans cette introduction.
\\ \\
La Th\'{e}orie des nombres \emph{hyperr\'{e}els} de A.Robinson est r\'{e}put\'{e}e difficile [2,3], sans doute parce qu'il s'agit de prouver par des arguments purement logiques l'existence d'un sur-espace de $\mathbb{R}$ qui soit un corps archim\'{e}dien, totalement ordonn\'{e} depuis l'infiniment petit jusqu'\`{a} l'infiniment grand.
\par
En particulier, on ne conna\^{i}t aucun de ces "hypernaturels" $\omega$,
on sait seulement qu'ils "existent" et d'une certaine mani\`{e}re,
 tout se passe comme lorsque l'on travaille avec des lettres,
 on n'en conna\^{i}t pas la valeur.
 \\
 Ici, l'ind\'{e}termination est tr\`{e}s
  limit\'{e}e puisqu'un seul \emph{nouveau} nombre va permettre d'exprimer tous les autres.
\\ \\
A l'oppos\'{e} de cette approche \emph{formaliste} de A.Robinson en Analyse Non Standard, on construit ici un sur-anneau de $\mathbb{R}$ qui n'est pas archim\'{e}dien, un nouvel ensemble de nombres entiers infiniment grands mais pas infinis au sens de G.Cantor et un sur-espace not\'{e} $\Omega$ de ces deux structures qui poss\`{e}de toutes les propri\'{e}t\'{e}s voulues.
\par
On ne saura ici d\'{e}nombrer que des ensembles de nombres tr\`{e}s particuliers et le plus grand espace num\'{e}rique ici consid\'{e}r\'{e}, $\overline\Omega$, n'est que le plus petit ensemble totalement ordonn\'{e} et continu qui contienne l'ensemble $\mathbb{R}$, diff\'{e}rent bien s\^{u}r de l'ensemble $\mathbb{R}$ lui-m\^{e}me.
\par
A titre de comparaison, l'ensemble des \emph{surr\'{e}els} \emph{construit} par Conway \`{a} base de coupures de Dedekind it\'{e}r\'{e}es \`{a} l'infini, conduit au plus gros sur-corps de $\mathbb{R}$ totalement ordonn\'{e}. Il contient \`{a} titre particulier, tous les ordinaux transfinis (l'addition y est commutative) ainsi que tous ses inverses.
\\ \\
On comprend alors la diff\'{e}rence de technicit\'{e} entre ces trois tentatives modernes de l\'{e}gitimer le r\^{e}ve Leibnizien.
\newpage
\emph{Toutes les branches des Math\'{e}matiques, y compris la M\'{e}canique classique et la Relativit\'{e} G\'{e}n\'{e}rale, lorsqu'elles recourent aux structures de nombres r\'{e}els ou entiers \emph{standard}, en particulier pour int\'{e}grer une \'{e}quation diff\'{e}rentielle, pourraient trouver b\'{e}n\'{e}fice \`{a} ces \'{e}largissements tr\`{e}s simples des ensembles $\mathbb{R}$ et $\mathbb{N}$ qui sont ici pour la premi\`{e}re fois propos\'{e}s.}
\\ \\
On trouvera bient\^{o}t dans ArXiv D.S. [5] une premi\`{e}re application importante de cette nouvelle Analyse Non Standard en Th\'{e}orie de la Gravitation (mais $\mathbb{R}_3=\mathbb{R}[X]/(X^3)$ remplace $\mathbb{R}[[X]]$ pour que les calculs soient plus simples).
\\ \\
Le contenu de cet article est le suivant :
\\ \\
Dans la premi\`{e}re partie, on \'{e}tudie les propri\'{e}t\'{e}s de l'extension de $\mathbb{R}$, not\'{e}e $\mathbb{R}_o$. On prolonge analytiquement les fonctions de classe $C^{\infty}$ de $\mathbb{R}$ dans $\mathbb{R}$ et l'on montre l'utilit\'{e} de ce prolongement pour le calcul des diff\'{e}rentielles. On d\'{e}finit \'{e}galement une Topologie pour l'\'{e}tude des propri\'{e}t\'{e}s diff\'{e}rentielles des fonctions de $\mathbb{R}_o$ dans $\mathbb{R}_o$.
\\
Dans la deuxi\`{e}me partie, on construit \emph{le} prolongement intrins\`{e}que de l'ensemble des nombres entiers standard, not\'{e} $\mathbb{N}[\Sigma]$, puis l'extension inductive $\aleph$ de $\mathbb{N}[\Sigma]$. La troisi\`{e}me partie est consacr\'{e}e au Th\'{e}or\`{e}me Fondamental et la quatri\`{e}me partie \`{a} l'\'{e}tude math\'{e}matique des oeuvres de I.Newton et G.W.Leibniz.
\\
Dans une cinqui\`{e}me et derni\`{e}re partie, on d\'{e}finit un corps archim\'{e}dien not\'{e} $\Omega$ qui est l'extension naturelle de l'ensemble $\mathbb{R}$ contenant \`{a} titre de sous-espaces $\mathbb{R}_o$ et $\aleph$ munis de leurs propri\'{e}t\'{e}s alg\'{e}briques et ordinales.

\newpage

\section{De nouveaux nombres "r\'{e}els".
\\ \\Premi\`{e}res propri\'{e}t\'{e}s de l'ensemble $\mathbb{R}_o$}

1.1. Propri\'{e}t\'{e}s alg\'{e}briques et ordinales de $\mathbb{R}_o=\mathbb{R}[[X]]$.
\\
1.2. Prolongement analytique d'une fonction de classe $C^\infty$ de $\mathbb{R}$ dans $\mathbb{R}$.
\\
1.3. "Continuit\'{e}", NS*-continuit\'{e} et "diff\'{e}rentiabilit\'{e}" multiple d'une fonction de $\mathbb{R}_o$ dans $\mathbb{R}_o$.
\\
1.4. "Diff\'{e}rentiabilit\'{e}" multiple du prolongement analytique d'une fonction $C^{\infty}$ standard.
\\
1.5. Deux types de diff\'{e}rentielles pour les prolongements analytiques des fonctions analytiques standard et leurs relations r\'{e}ciproques.
\\
1.6. Sur le contexte historique de cette recherche.
\\
1.7. Le co\^{u}t d'une r\'{e}volution non standard.

\newpage
\subsection{Propri\'{e}t\'{e}s alg\'{e}briques et ordinales de $\mathbb{R}_o=\mathbb{R}[[X]]$}
L'ensemble des s\'{e}ries formelles \`{a} une ind\'{e}termin\'{e}e
 \begin{center}
$\mathbb{R}[[X]]=\lbrace{x=\sum\limits_{k\geq{0}}a_kX^k}/a_k\in{\mathbb{R}}\rbrace$
\end{center}
 muni des lois usuelles est bien s\^{u}r un espace vectoriel de dimension infinie d\'{e}nombrable sur $\mathbb{R}$ et un anneau commutatif unitaire.

\begin{prop} $(\mathbb{R}[[X]],+,\times)$ est un anneau int\`{e}gre.
\end{prop}
Preuve : soit $x=\sum\limits_{k\geq{ord(x)}}a_kX^k$, $ord(x)$ est le plus petit entier $k$ tel que $a_k\neq{0}$. Par convention, $ord(0)=+\infty$. On montre que $ord(x\times{y})=ord(x)+ord(y)$. Si $x\times{y}=0$ et $x\neq{0}$, $y\neq{0}$, contradiction. \qed

\begin{rem} Les coefficients $a_k$ peuvent \^{e}tre en nombre infini non nuls mais les puissances $k$ de $X$ dans $X^k$ sont toutes finies.
\end{rem}

\begin{defn} On d\'{e}finit un ordre total sur $\mathbb{R}[[X]]$. C'est l'ordre lexicographique :
\begin{center}
$x=y$ ssi $\forall{n}\in{\mathbb{N}}$ $a_n=b_n$
\\Sinon, soit $p=ord(y-x)$, i.e. $(\forall{n<p})$ $a_n=b_n$ et $a_p\neq{b_p}$.
\\Alors $x<y$ ssi $a_p<b_p$. On a $x\leq{y}$ ssi $(x=y)$ ou $(x<y)$.
\end{center}
\end{defn}

\begin{prop} $(\mathbb{R}[[X]],+,\cdot,\times,\leq)$ est une alg\`{e}bre totalement ordonn\'{e}e.
\end{prop}
Preuve : on d\'{e}montre seulement que $\leq$ est compatible avec la multiplication interne (ce n'est pas le cas dans $\mathbb{C}$ puisque, si $i$ \'{e}tait positif, on aurait $i^2=-1$ positif). Si $x=\sum\limits_{k\geq{p}}a_kX^k>0$ et $y=\sum\limits_{l\geq{q}}b_lX^l>0$, avec $p=ord(x)$ et $q=ord(y)$, on a $a_p>0$ et $b_q>0$ donc $a_pb_q>0$ et $x\times{y}>0$.\qed

\begin{rem} 1) On note la structure totalement ordonn\'{e}e de $\mathbb{R}[[X]]$ par $\mathbb{R}_o$, en souvenir de I.Newton ([6] et [7, p.261]) et l'on consid\`{e}re les \'{e}l\'{e}ments de la structure d'alg\`{e}bre totalement ordonn\'{e}e $(\mathbb{R}_o,+,\cdot,\times,\leq)$ comme aussi "r\'{e}els" que les \'{e}l\'{e}ments de la structure standard \hbox{$(\mathbb{R},+,\times,\leq)$}.
\par 2) On ne peut recourir \`{a} la notation $\mathbb{R}[o]$ ou $\mathbb{R}(o)$ car elles sont utilis\'{e}es classiquement [8, p.182] pour d\'{e}signer le plus petit sur-anneau et sur-corps de $\mathbb{R}$ qui contienne l'\'{e}l\'{e}ment transcendant nouveau $o$, c'est l'ensemble des polyn\^{o}mes de degr\'{e} fini et l'ensemble des fonctions rationnelles du nombre $o$.
\end{rem}
$\mathbb{R}_o$ n'est pas un corps. Il n'est pas non plus archim\'{e}dien car, selon l'ordre lexicographique pr\'{e}c\'{e}dent :
\begin{center}
$(\forall {k\in\mathbb{N}})$ $ko<1$, not\'{e} $0<o\ll{1}$
\end{center}
On dit que $o$ est infiniment petit, ou infinit\'{e}simal. On note aussi
\begin{center}
$x\ll{y}$ ssi $(\forall{k}\in\mathbb{N})$ $k\vert{x}\vert<\vert{y}\vert$
\end{center}
$x_S=t=a_0$ est la partie standard du nombre r\'{e}el $x$.
\\$u=x-x_S$ est sa partie infinit\'{e}simale.
\\$]x[=\lbrace{y\in\mathbb{R}_{o}}/y_S=x_S\rbrace$ est la coupure infinit\'{e}simale de $\mathbb{R}_{o}$ en $x$ et $[0[$ la demi-coupure \`{a} droite en $0$, c'est l'ensemble des nombres infinit\'{e}simaux positifs de $\mathbb{R}_{o}$.
\\
Tout \'{e}l\'{e}ment de $\mathbb{R}_{o}$ s'\'{e}crit donc de mani\`{e}re unique $x=t+u$, avec $t\in{\mathbb{R} }$ et $\vert{u}\vert\ll{1}$, $u\in{]0[}$ et $a_k\cdot{o^k}$ est le moment d'ordre $k$ de $u=\sum\limits_{k\geq{1}}a_k\cdot{o^k}$. $T_N(x)=\sum\limits_{0\leq{k}\leq{N}}a_k\cdot{o^k}$ est la troncature \`{a} l'ordre $N$ du nombre r\'{e}el non standard $x$.

\subsection{Prolongement analytique d'une fonction de classe $C^\infty$ de $\mathbb{R}$ dans $\mathbb{R}$}
\emph{$\mathbb{R}_o$ \'{e}tant muni de sa Topologie d'ordre, la plupart des fonctions de $\mathbb{R}$ dans $\mathbb{R}_o$ ne sont pas continues en $t$, a fortiori d\'{e}rivables (on montre facilement qu'une fonction continue de $\mathbb{R}$ dans $\mathbb{R}_o$ a des troncatures constantes \`{a} n'importe quel ordre sur des ouverts standard de $t_0$ de plus en plus petits). C'est pourquoi on prolonge de mani\`{e}re analytique seulement les fonctions $C^{\infty}$ de $\mathbb{R}$ dans $\mathbb{R}$, dans un premier temps.
\\
 En 1.3 et 3.2, on \'{e}tudie les propri\'{e}t\'{e}s de classes de fonctions de $\mathbb{R}_o$ dans $\mathbb{R}_o$ (NS*-continues, $p$ fois "diff\'{e}rentiables" et r\'{e}guli\`{e}res ou newtoniennes). La sous-section 1.5 est \'{e}crite pour les prolongements analytiques des fonctions analytiques standard mais elle se g\'{e}n\'{e}ralise sans changement aux fonctions r\'{e}guli\`{e}res qui sont d\'{e}finies en 3.2 puisque ces fonctions ont presque les m\^{e}mes propri\'{e}t\'{e}s que les prolongements analytiques.}
\\ \\
R.Godement [7, p.264-7] \'{e}tudie rapidement le prolongement analytique d'une fonction de $\mathbb{R}$ dans $\mathbb{R}$, de classe $C^{I}$, dans l'anneau \emph{non int\`{e}gre} $\mathbb{R}[X]/(X^{I})$, $I=2$ ou $3$. Cf. aussi [5].
\par
Ici, tout se passe de la m\^{e}me mani\`{e}re sauf que la S\'{e}rie de Taylor peut avoir un nombre infini d\'{e}nombrable de termes non nuls.
\\ \\
Soit une fonction $f$ de $\mathbb{R}$ dans $\mathbb{R}$ de classe $C^{\infty}$. M\^{e}me si sa S\'{e}rie de Taylor a un rayon de convergence $R=0$ en $t$, celle-ci converge \`{a} l'int\'{e}rieur de la coupure infinit\'{e}simale de $t$. En effet, le nombre $\sum\limits_{k\geq{0}}\frac{1}{k!}f^{(k)}(t)\cdot{u^k}$ est \emph{un} \'{e}l\'{e}ment bien d\'{e}fini de $\mathbb{R}_{o}$.

\begin{defn} Le prolongement analytique d'une fonction $f$ de classe $C^\infty$, est la fonction not\'{e}e $\bar{f}$ de $\mathbb{R}_{o}$ dans $\mathbb{R}_{o}$, telle que
\begin{center}
$\bar{f}(t+u)=\sum\limits_{k\geq{0}}\frac{1}{k!}f^{(k)}(t)\cdot{u^k}$.
\end{center}
\end{defn}

\begin{rem} On peut calculer dans $\mathbb{R}_{o}$ tous les coefficients des puissances finies de $o$, avec $u=\sum\limits_{j\geq{1}}a_j\cdot{o^j}$. On consid\`{e}re :
\begin{center}
$\bar{f}_N(t+u)=\sum\limits_{0\leq{k}\leq{N}}\frac{1}{k!}f^{(k)}(t)
[\sum\limits_{1\leq{j}\leq{N}}a_j\cdot{o^j}]^k\cdot{o^{i}}$
\end{center}
qui a la m\^{e}me partie standard et les $N$ premiers moments infinit\'{e}simaux que $\bar{f}(t+u)$.
\end{rem}

\begin{prop} $f:\mathbb{R}^+\longrightarrow{\mathbb{R}^+}$, $t\longmapsto{t^\alpha}$ avec $\alpha>0$, $\mathbb{R}^{+*}\longrightarrow{\mathbb{R}^{+*}}$ si $\alpha<0$. On prolonge analytiquement $f$ dans $\mathbb{R}_{o}^+$ si $\alpha>0$, dans $\mathbb{R}_{o}^+\backslash{[0[}$ si $\alpha<0$
et $\bar{f}(x)=t^\alpha\sum\limits_{k\geq{0}}C_\alpha^k\cdot{(\frac{u}{t})^k}$ avec $C_\alpha^k=\frac{\alpha{(\alpha-1)}\dots(\alpha-k+1)}{k!}$ m\^{e}me si $\alpha$ n'est pas un entier.
\end{prop}

\begin{cor} Tout nombre non infinit\'{e}simal est inversible et $\frac{1}{t+u}=\frac{1}{t}\sum\limits_{k\geq{0}}(-1)^k(\frac{u}{t})^k$.
\end{cor}

\subsection{"Continuit\'{e}", NS*-continuit\'{e} et "diff\'{e}rentiabilit\'{e}" multiple d'une fonction de $\mathbb{R}_o$ dans $\mathbb{R}_o$}
\emph{On \emph{ne} peut d\'{e}finir pour une fonction \`{a} valeurs r\'{e}elles non standard, une notion classique de d\'{e}riv\'{e}e car $u$ n'est pas inversible. $\mathbb{R}_o$ n'est pas un Espace de Banach et l'on ne peut profiter de la g\'{e}n\'{e}ralisation introduite par J.Dieudonn\'{e} [9] dans le Calcul Infinit\'{e}simal.}
\par
\emph{On pourrait d\'{e}finir une simple op\'{e}ration de d\'{e}rivation dans l'ensemble des fonctions $\bar{f}$ par l'\'{e}galit\'{e} $\bar{f}$ $'=\bar{f'}$. On pr\'{e}f\`{e}re d\'{e}finir en 1.4 une notion g\'{e}n\'{e}rale de diff\'{e}rentiabilit\'{e} qui soit compatible avec la Topologie d'ordre de $\mathbb{R}_o$. D'abord, on d\'{e}finit une notion sp\'{e}cifique \`{a} l'Analyse Non Standard [3], la NS*-continuit\'{e}.}

\begin{lem} Une fonction $F$ de $\mathbb{R}_o$ dans $\mathbb{R}_o$ est NS*-continue si elle v\'{e}rifie l'une des deux propri\'{e}t\'{e}s \'{e}quivalentes :
\par 1) Pour tous les $x_1,x_2$ de $\mathbb{R}_o$, $|x_2-x_1|\ll{|F(x_2)-F(x_1)|}$ n'est jamais v\'{e}rifi\'{e}.
\par 2) $(\forall{k\in{\mathbb{N}^*}})(\forall{x\in{\mathbb{R}_o}})$ $F(]x[_k)\subseteq{]F(x)[_k}$. $]y[_k$ est l'ensemble des nombres proches de $y$ \`{a} un nombre infinit\'{e}simal pr\`{e}s, d'ordre au moins $k$ (c'est la coupure en $y$ d'ordre $k$).
\end{lem}

\begin{defn} $F$ est "continue" en $x_0$ ssi :
\\$(\forall{\varepsilon\in{\mathbb{R}_o^{+*}}})(\exists{\eta\in{\mathbb{R}_o^{+*}}})(\forall{x\in{\mathbb{R}_o}})$
$|x-x_0|<\eta\Rightarrow{|F(x)-F(x_0)|<\varepsilon}$.
\end{defn}

\begin{rem} L'ast\'{e}risque signifie que la NS*-continuit\'{e} n'est pas celle de l'Analyse Non Standard de Robinson. Les guillemets signalent qu'il ne s'agit pas des notions de l'Analyse standard.
\end{rem}

\begin{prop} Si $F$ est NS*-continue, alors $F$ est "continue" partout.
\end{prop}
Preuve : soit $\varepsilon=o^n$. On prend $\eta=o^{n+1}$, alors $ord(F(x)-F(x_0))\geq{ord(x-x_0)}\geq{n+1}$ donc $|F(x)-F(x_0)|<o^n$.

\begin{rem} La r\'{e}ciproque est fausse. Par exemple $F$ de $\mathbb{R}_o$ dans $\mathbb{R}_o$ telle que $x=t+a_1\cdot{o}+a_2\cdot{o^2}+\dots{\mapsto{F(x)=a_1+a_2\cdot{o}+\dots}}$ est partout "continue" (prendre $\eta=\varepsilon\times{o}$) mais pas NS*-continue.
\end{rem}

\begin{defn} Une fonction $F$ de $\mathbb{R}_o$ dans $\mathbb{R}_o$ est $p$ fois "diff\'{e}rentiable" en $x_0\in{\mathbb{R}_o}$ s'il existe $p+1$ nombres $H_i\in{\mathbb{R}_o}$, $i=0,\dots{p}$ tels que :
\\$(\forall{\varepsilon}\in\mathbb{R}_{o}^{+*})
(\exists{\eta}\in\mathbb{R}_{o}^{+*})
(\forall{x}\in\mathbb{R}_{o})\vert{x-x_0}\vert{<\eta}
\Rightarrow{\vert{F(x)}}-H_0-H_1\times{(x-x_0)-\dots{H_p\times{(x-x_0)^p}}}
\vert<\varepsilon\vert{x-x_0}\vert^p$.
\end{defn}

\begin{prop} 1) Si $F$ est $p+1$ fois "diff\'{e}rentiable" en $x_0$ alors $F$ est $p$ fois "diff\'{e}rentiable".
\par 2) Si $F$ est "diff\'{e}rentiable" en $x_0$ alors $F$ est NS*-continue \hbox{en $x_0$} (fixer $x_2=x_0$ dans la d\'{e}finition 1 ou $x=x_0$ dans la d\'{e}finition 2).
\par 3) S'ils existent, ces nombres $H_i$ sont uniques et d\'{e}pendent en g\'{e}n\'{e}ral de $x_0$.
\end{prop}

\subsection{"Diff\'{e}rentiabilit\'{e}" du prolongement analytique d'une fonction $C^{\infty}$ standard}

\begin{rem} Cette notion non standard* de "diff\'{e}rentiabilit\'{e}" est locale, infinit\'{e}simale m\^{e}me. Puisqu'une fonction $\bar{f}$ est d\'{e}finie par une s\'{e}rie formelle, on d\'{e}montre qu'elle est "diff\'{e}rentiable" partout.
\end{rem}

\begin{prop} Tout prolongement analytique d'une fonction num\'{e}rique infiniment d\'{e}rivable, est "diff\'{e}rentiable" en $t\in{\mathbb{R}}$.
\end{prop}
Preuve : on prend $\varepsilon$ et $\eta$ infinit\'{e}simaux et on cherche $H$ tel que
\\$\vert{\bar{f}(t+u)-\bar{f}(t)-H\times{u}}\vert<\varepsilon\times{\vert{u}\vert}$ pour $\vert{u}\vert<\eta$.
\\ On prend bien s\^{u}r $H=f'(t)$ et $\vert{\sum_{k\geq{2}}}\frac{f^{(k)}(t)}{k!}\cdot{u^k}\vert<b_2\cdot{u^2}$ pour tout r\'{e}el standard $b_2$ strictement positif tel que $\vert{\frac{1}{2}f''(t)}\vert<b_2$, du fait de l'ordre lexicographique.
\\ On a alors $(\forall{\varepsilon}\in{[0[})(\exists{\eta}\in{[0[})(\forall{u}\in{]0[})\vert{u}\vert<\eta
\Rightarrow{b_2\cdot{u^2}}<\varepsilon{\vert{u}}\vert$.
\\Il suffit de prendre $\eta=\frac{\varepsilon}{b_2}$ \qed

\begin{defn} On dit alors que $\bar{f}$ est "d\'{e}rivable" en $t$. De mani\`{e}re g\'{e}n\'{e}rale, on note $H=F'(t)$ et, de mani\`{e}re particuli\`{e}re aux prolongements analytiques, on a $\bar{f}'=\bar{f'}$.
\end{defn}

\begin{prop} Tout prolongement analytique d'une fonction num\'{e}rique infiniment d\'{e}rivable est "d\'{e}rivable" en $x\in{\mathbb{R}_o}$, \hbox{$x=t+u$}.
\end{prop}
Preuve : on prend bien s\^{u}r $H=\sum\limits_{k\geq{0}}\frac{f^{(k+1)}(t)}{k!}\cdot{u^k}$ et
$\vert{\bar{f}(x+v)}-\bar{f}(x)-H\times{v}\vert=\vert{\sum\limits_{k\geq{2}}}\frac{f^{(k)}}{k!}
\cdot{[(u+v)^k-u^k-ku^{k-1}v]}<b_2\cdot{u^2}$ comme pr\'{e}c\'{e}demment. \qed

\begin{rem} 1) On a $\bar{f}'(x)=\sum\limits_{k\geq{0}}\frac{f^{(k+1)}(t)}{k!}\cdot{u^k}$. On remarque qu'on ne d\'{e}rive que par rapport \`{a} la partie standard $t$ de $x$. On dira donc que cette op\'{e}ration de "d\'{e}rivation" fournit la d\'{e}riv\'{e}e partielle de $\bar{f}$ par rapport \`{a} la partie standard de $x$.
\par 2) On d\'{e}montre de m\^{e}me que $\bar{f}$ est infiniment "diff\'{e}rentiable" en $x$ et l'on d\'{e}finit par r\'{e}currence les "d\'{e}riv\'{e}es" successives de la fonction $\bar{f}$. On retrouve la formule bien connue de l'Analyse standard :
$$\bar{f}^{(q)}(x)=\sum\limits_{k\geq{0}}\frac{f^{(k+q)}(t)}{k!}\cdot{u^k}$$ et $\bar{f}^{(q)}=\bar{f^{(q)}}$.
\end{rem}

\begin{prop} On a $\bar{f}(x+v)=\sum\limits_{q\geq{0}}\frac{1}{q!}\bar{f}^{(q)}(x)\times{v^q}$ pour $x\in\mathbb{R}_{o}$ et $\vert{v}\vert\ll{1}$.
\end{prop}
Preuve : on prouve l'\'{e}galit\'{e} jusqu'au moment infinit\'{e}simal d'ordre $N$. On a $\bar{f}_N(x+v)=\sum\limits_{0\leq{k}\leq{N}}\frac{1}{k!}f^{(k)}(t)\cdot{(u+v)^k}$.
\\
Classiquement,
$\bar{f}_N(x+v)=\sum\limits_{0\leq{k}\leq{N}}\sum\limits_{0\leq{q}\leq{k}}f^{(k)}(t)\cdot{(\frac{u^{k-q}}{(k-q)!}
\times{\frac{v^q}{q!}})}$ =
\\
$\sum\limits_{0\leq{q}\leq{N}}\frac{1}{q!}[\sum\limits_{q\leq{k}\leq{N}}\frac{1}{(k-q)!}f^{(k)}(t)\cdot{u^{k-q}}]
\times{v^q}$ qui a le m\^{e}me d\'{e}but que \hbox{$\sum\limits_{q\leq{N}}\frac{1}{q!}\bar{f}^{(q)}(x)\times{v^q}$}. \qed

\subsection{Deux types de diff\'{e}rentielles pour les prolongements analytiques des fonctions analytiques standard et leurs relations r\'{e}ciproques}
\begin{lem} (Th\'{e}or\`{e}me des diff\'{e}rences finies [10, p.204])
\\
Soit $(u_n)_{n\in{\mathbb{N}}}$, une suite r\'{e}elle. $\Delta{u_n}=u_{n+1}-u_n$ et $\Delta^{p+1}u_n=\Delta^pu_{n+1}-\Delta^pu_n$. On a $\Delta^pu_n=\sum\limits_{k=0}^p(-1)^{p-k}C_p^ku_{n+k}$.
\end{lem}
Preuve : par r\'{e}currence \`{a} partir de $\Delta^2u_n=(u_{n+2}-u_{n+1})-(u_{n+1}-u_n)=u_{n+2}-2u_{n+1}+u_n$. \qed
\\ \\
Il en est de m\^{e}me du prolongement analytique d'une fonction $f$ de $\mathbb{R}$ dans $\mathbb{R}$, de classe $C^\infty$, \`{a} l'int\'{e}rieur d'une coupure infinit\'{e}simale $]t[$ pour l'op\'{e}rateur de diff\'{e}rentiation $D^p$ d\'{e}fini par r\'{e}currence comme $\Delta^p$.

\begin{defn} On appelle $p$-i\`{e}me diff\'{e}rentielle de $f$ et on note $D^p\bar{f}$ la fonction de $\mathbb{R}_{o}$ dans $]0[$ d\'{e}finie par r\'{e}currence par
\begin{center}
$D\bar{f}(x)=\bar{f}(x+o)-\bar{f}(x)$ et $D^{p+1}\bar{f}(x)=D^p\bar{f}(x+o)-D^p\bar{f}(x)$.
\end{center}
\end{defn}

On d\'{e}montre imm\'{e}diatement que $(D\bar{f})'=D\bar{f'}$ et, dans le cas g\'{e}n\'{e}ral $(DF)'=DF'$.
\begin{prop} On a :
\\$D^p\bar{f}(x)=\sum\limits_{k=0}^p(-1)^{p-k}C_p^k\cdot{\bar{f}(x+ko)}=\sum\limits_{n\geq{0}}\frac{X_p^n(u)}{n!}f^{(n)}(t)$ avec
\\$X_p^n(u)=\sum\limits_{k=0}^p(-1)^{p-k}C_p^k\cdot{(u+ko)^n}$ et $t=x_S$, $u=x-x_S$.
\end{prop}
Preuve : La d\'{e}monstration est la m\^{e}me que celle du \textbf{Lemme 1.14}. \qed

\begin{cor} $X_p^n(u)=0$ si $n<p$ et $X_p^p(u)=p!\cdot{o^p}$.
\end{cor}
Preuve : on \'{e}crit la relation pr\'{e}c\'{e}dente pour un polyn\^{o}me $\bar{f}$ quelconque de degr\'{e} $p-1$. $D^p\bar{f}(x)=0$ car le degr\'{e} du polyn\^{o}me diminue d'une unit\'{e} \`{a} chaque \emph{diff\'{e}renciation}.
\\
On a aussi $f^{(n)}(t)=0$ si $n\geq{p}$ et donc aucune condition sur $X_p^n(u)$ dans ce cas. Par contre, si $n<p$, il faut que les coefficients $X_p^n(u)$ soient tous nuls pour que la relation $D^p\bar{f}(x)=0$ soit v\'{e}rifi\'{e}e quelle que soit la fonction $\bar{f}$.
\\
Enfin, $D^{p-1}\bar{f}(x)=a_{p-1}(p-1)!o^{p-1}$ et $f^{(p-1)}(t)=a_{p-1}(p-1)!$ o\`{u} $a_{p-1}$ est le coefficient du mon\^{o}me dominant de $f$. Donc :
\\$\sum\limits_{n\geq{0}}\frac{X_{p-1}^n(u)}{n!}f^{(n)}(t)=0+\frac{X_{p-1}^{p-1}(u)}{(p-1)!}
f^{(p-1)}(t)+0=D^{p-1}\bar{f}(x)$ et
\\ $X_{p-1}^{p-1}(u)=(p-1)!\cdot{o^{p-1}}$
 \qed

\begin{defn} On appelle diff\'{e}rentielle d'ordre $n$ not\'{e}e $d^n\bar{f}$, la fonction de $\mathbb{R}_{o}$ dans $]0[$ d\'{e}finie depuis G.W.Leibniz par
\begin{center}
$d^n\bar{f}(x)=\bar{f}^{(n)}(x)\times{o^n}$ ($o=dx$).
\end{center}
\end{defn}

\begin{rem} 1) On n'utilise pas la notation de Leibniz
$$\frac{d^n\bar{f}}{dx^n}=\bar{f}^{(n)}$$
car $o$ n'est pas inversible.
\par 2) D'apr\`{e}s la \textbf{D\'{e}finition 1.24}, $D\bar{f}(x)=\sum_{n\geq{0}}\frac{d^n\bar{f}(x)}{n!}$.
\end{rem}

\begin{prop} On a $D^p\bar{f}(x)=\sum\limits_{n\geq{p}}X_p^n\times{\frac{d^n\bar{f}(x)}{n!}}$ avec
$X_p^n=\sum\limits_{k=0}^p(-1)^{p-k}C_p^kk^n$.
\end{prop}
Cette relation donne les diff\'{e}rentielles $p$-i\`{e}mes en fonction des diff\'{e}rentielles d'ordre $n$ avec $n\leq{p}$. On donne plus loin les relations r\'{e}ciproques.

\begin{rem} En particulier :
 \\$D\bar{f}(x)=d\bar{f}(x)+\frac{1}{2}d^2\bar{f}(x)+\frac{1}{6}\cdot{d^3\bar{f}(x)}+
 \frac{1}{24}\cdot{d^4\bar{f}(x)}+\dots$, \\$D^2\bar{f}(x)=d^2\bar{f}(x)+d^3\bar{f}(x)+\frac{7}{12}\cdot{d^4\bar{f}(x)}+\dots$
\\$D^3\bar{f}(x)=d^3\bar{f}(x)+\frac{3}{2}\cdot{d^4\bar{f}(x)}+\dots$.
\end{rem}

\begin{prop} On a $\frac{d^n\bar{f}(x)}{n!}=
\sum\limits_{p\geq{n}}(-1)^{p-n}K_{p-1}^{p-n}\cdot{\frac{D^p\bar{f}(x)}{p!}}$ avec $x\in{\mathbb{R}_o}$. $K_{p-1}^{p-n}$ est la somme des $C_{p-1}^{p-n}=C_{p-1}^{n-1}$ produits possibles de $p-n$ facteurs pris parmi les $p-1$ premiers entiers non nuls.
\end{prop}
Preuve : On a $\bar{f}_N(x+ko)=
\bar{f}(x)+\sum\limits_{n=1}^N\frac{d^n\bar{f}(x)}{n!}k^n$
o\`{u} N et k sont deux entiers
 finis suffisamment grands.
  On d\'{e}montre par r\'{e}currence finie que $\bar{f}_N(x+ko)=
  \bar{f}(x)+\sum\limits_{p=1}^NC_k^p\cdot{D^p\bar{f}(x)}$.
\\
Si $p\geq{2}$, $C_k^p=\frac{k(k-1)
\dots{(k-(p-1))}}{p!}=\frac{1}{p!}
[\sum\limits_{n=1}^{p}(-1)^{p-n}K_{p-1}^{p-n}k^n]$
et
\\
$\bar{f}_N(x+ko)=\bar{f}(x)+k\cdot{D\bar{f}(x)}+\sum\limits_{2\leq{p}\leq{N}}[\sum\limits_{1\leq{n}\leq{p}}
(-1)^{p-n}K_{p-1}^{p-n}k^n]\cdot{\frac{D^p\bar{f}(x)}{p!}}=\bar{f}(x)+k\cdot{D\bar{f}(x)}+\sum\limits_{1\leq{n}\leq{N}}[\sum\limits_{n\leq{p}\leq{N}}^{p\geq{2}}
(-1)^{p-n}K_{p-1}^{p-n}\cdot{\frac{D^p\bar{f}(x)}{p!}}]k^n$.
\\
Le coefficient de $k$ (pour $n=1$) est $D\bar{f}(x)-\frac{D^2\bar{f}(x)}{2}+\frac{D^3\bar{f}(x)}{3}-\dots$ car $K_{p-1}^{p-1}=(p-1)!$. Il est \'{e}gal \`{a} $d\bar{f}(x)$ d'apr\`{e}s la premi\`{e}re \'{e}galit\'{e}.
\\
Les deux coefficients de $k^n$ (pour $n\geq{2}$) sont \'{e}gaux et donc
\begin{center}
$\frac{d^n\bar{f}_N(x)}{n!}=\sum\limits_{n\leq{p}\leq{N}}(-1)^{p-n}K_{p-1}^{p-n}\cdot{\frac{D^p\bar{f}(x)}{p!}}$
\end{center}
jusqu\`{a} l'ordre $N$. On prolonge les \'{e}galit\'{e}s jusqu'\`{a} l'infini (pour toute valeur finie $N$)
\begin{center}
$\frac{d^n\bar{f}(x)}{n!}=\sum\limits_{n\leq{p}}(-1)^{p-n}K_{p-1}^{p-n}\cdot{\frac{D^p\bar{f}(x)}{p!}}$.
\end{center} \qed

\begin{rem} On trouve en particulier :
\\$d\bar{f}(x)=D\bar{f}(x)-\frac{D^2\bar{f}(x)}{2}+\frac{D^3\bar{f}(x)}{3}
-\frac{D^4\bar{f}(x)}{4}+\dots$
\\$d^2\bar{f}(x)=D^2\bar{f}(x)-D^3\bar{f}(x)+\frac{11}{12}\cdot{D^4\bar{f}(x)}+\dots$
\\$d^3\bar{f}(x)=D^3\bar{f}(x)-\frac{3}{2}\cdot{D^4\bar{f}(x)}+\dots$.
\end{rem}

\subsection{Sur le contexte historique de cette recherche}
On donne ici deux citations, l'une est de N.Bourbaki, l'autre est de R.Godement :
\begin{quotation}
"il faut bien reconna\^{i}tre que la notation leibnizienne de diff\'{e}rentielle n'a \`{a} vrai dire aucun sens ; au d\'{e}but du XIX\`{e}me si\`{e}cle, elle tomba dans un discr\'{e}dit dont elle ne s'est relev\'{e}e que peu \`{a} peu ; et, si l'emploi des diff\'{e}rentielles premi\`{e}res a fini par \^{e}tre compl\`{e}tement l\'{e}gitim\'{e}, les diff\'{e}rentielles d'ordre sup\'{e}rieure, d'un usage pourtant si commode, n'ont pas encore \'{e}t\'{e} vraiment r\'{e}habilit\'{e}es jusqu'\`{a} ce jour" [10, p.216] et [11].
\\ \\
"Ces notions qui reposent sur des "infiniment petits" que personne n'a jamais pu d\'{e}finir, ont fait inutilement cogiter et divaguer beaucoup trop de gens pour qu'on leur attribue maintenant un autre r\^{o}le que celui d'une explication historique de la notation diff\'{e}rentielle" [7, p.260].
\end{quotation}
Ce probl\`{e}me ancien a \'{e}t\'{e} je crois, ici r\'{e}solu au b\'{e}n\'{e}fice de la rigueur "et" de la compr\'{e}hension.
\\ \\
L'impossibilit\'{e} \`{a} un moment de l'Histoire de satisfaire \`{a} la fois aux exigences de la rigueur \textit{et} \`{a} une certaine "\'{e}vidence" des r\'{e}sultats semble \^{e}tre l'un des leitmotiv des \underline{El\'{e}ments d'Histoire des} \underline{ Math\'{e}matiques} de N.Bourbaki [10]. Ce conflit a bien s\^{u}r toujours \'{e}t\'{e} r\'{e}solu au b\'{e}n\'{e}fice de la rigueur mais on peut aujourd'hui reprendre cette question et tenter de concilier rigueur et \'{e}vidence dans les d\'{e}monstrations.
\\ \\ Quelques exemples :

\par On a justifi\'{e} en 1.1 et 1.2 \underline{La m\'{e}thode des fluxions} [6] de I.Newton en d\'{e}finissant une alg\`{e}bre totalement ordonn\'{e}e $\mathbb{R}_o$ avec laquelle on peut prolonger les fonctions $C^{\infty}$ par des s\'{e}ries enti\`{e}res. Les propri\'{e}t\'{e}s d'espace vectoriel de $\mathbb{R}_o$ permettent de justifier les \'{e}galisations faites par I.Newton entre deux expressions \`{a} l'ordre $n$.
  \par
  En 1.3, on a \'{e}tabli le lien entre Diff\'{e}rences finies et deux sortes de Diff\'{e}rentielles ($D^p\bar{f}$ et $d^n\bar{f}$) que G.W.Leibniz a seulement intuitivement per\c{c}u.
\\ \\
D'apr\`{e}s N.Bourbaki en effet :
  \begin{quotation} "il se tient tr\`{e}s pr\`{e}s du calcul des diff\'{e}rences finies dont son calcul diff\'{e}rentiel se d\'{e}duit pas un passage \`{a} la limite que bien entendu, il serait en peine de justifier rigoureusement ; et par la suite il insiste volontiers sur le fait que les principes s'appliquent indiff\'{e}remment \`{a} l'un et \`{a} l'autre" [10, p.208].
\end{quotation}

On peut appeler \emph{Histonique} (l'Histonique est \`{a} l'Histoire des Sciences ce que la Bionique est \`{a} la Biologie) cette nouvelle mani\`{e}re de faire des Math\'{e}matiques (et de la Physique, cf. 4.4 et [5]) au del\`{a} des probl\'{e}matiques historiques. Il s'agit de reprendre les oeuvres des grands math\'{e}maticiens des si\`{e}cles pass\'{e}s d'un point de vue "moderne".
\par
Dans le m\^{e}me esprit, on montre en 2.5 qu'une nouvelle approche "cantorienne" des nombres entiers rejoint l'intuition newtonienne et en 3. que Leibniz a eu raison de percevoir une analogie entre Int\'{e}grale et Somme, Diff\'{e}rentielle et Diff\'{e}rence.

\newpage
\subsection{Le co\^{u}t d'une r\'{e}volution non standard}
E.Benoit, sp\'{e}cialiste fran\c{c}ais de l'Analyse Non Standard de A.Robinson, r\'{e}sume ainsi la situation actuelle :
\begin{quotation}
"L'immensit\'{e} du savoir math\'{e}matique \'{e}nonc\'{e} dans le formalisme classique de Cauchy-Bolzano-Weierstrass-Cantor-... rend important le co\^{u}t d'une r\'{e}volution non standard qu'on pourrait esp\'{e}rer. C'est la raison fondamentale qui fait que cette nouvelle th\'{e}orie n'obtient pas dans la communaut\'{e} math\'{e}matique la place \`{a} laquelle elle pourrait pr\'{e}tendre".
 \\ in page personnelle in http://perso.univ-lr.fr/ebenoit/$\sharp{ans}$.
\end{quotation}
La situation est un peu diff\'{e}rente pour cette Nouvelle Analyse Non Standard (que l'on pourrait noter ANS*) car elle appartient au cadre classique. C'est en effet en utilisant les notions standard de la Topologie G\'{e}n\'{e}rale que l'on a d\'{e}montr\'{e} les propri\'{e}t\'{e}s de $\mathbb{R}_o$.
\\ \\
On fait ici un rapide bilan du changement d'habitudes mentales qui serait n\'{e}cessaire pour adopter en Math\'{e}matiques, l'ensemble $\mathbb{R}_o$.

\paragraph{1} La propri\'{e}t\'{e} essentielle de $\mathbb{R}$ qui manque \`{a} $\mathbb{R}_o$ est l'existence d'une borne sup\'{e}rieure \`{a} tout sous-ensemble major\'{e} (Th\'{e}or\`{e}me de Bolzano). Ainsi, $[0[$ est major\'{e} par tous les nombres standard strictement positifs. On pourra y suppl\'{e}er partiellement \`{a} partir de 5.4 avec l'\'{e}tude de la continuit\'{e}.
\\
Cette perte entraine celle des Th\'{e}or\`{e}mes de Borel-Lebesgue et de Bolzano-Weierstrass mais surtout la perte du Th\'{e}or\`{e}me des valeurs interm\'{e}diaires qui est pourtant tr\`{e}s intuitif.
\par
Par exemple, la simple fonction $x\mapsto{x^2}$ envoie l'intervalle $[0,1]$ sur la r\'{e}union disjointe $[0[[\cup{[0,1]\setminus{[0[}}$. On montre qu'un prolongement analytique peut ne pas v\'{e}rifier la propri\'{e}t\'{e} des valeurs interm\'{e}diaires seulement aux points $x\in{\mathbb{R}_o}$ tels que $f'(x_S)=0$.

\begin{prop} $f:\mathbb{R}\rightarrow{\mathbb{R}}$, $C^{\infty}$ et $\bar{f}:\mathbb{R}_o\rightarrow{\mathbb{R}_o}$. On a $\bar{f}(x_1)=y_1$ et $\bar{f}(x_2)=y_2$ avec $x_1^S<x_2^S$ et $y_1^S<y_2^S$. Soit $y\in{]y_1,y_2[}$, si $f$ est strictement monotone sur $]x_1,x_2[$, il existe au moins un $x\in{]x_1,x_2[}$ tel que $\bar{f}(x)=y$.
\end{prop}
Preuve : on a $f(x_1^S)=y_1^S$ et $f(x_2^S)=y_2^S$. La fonction $f$ v\'{e}rifie la propri\'{e}t\'{e} des valeurs interm\'{e}diaires, donc il existe $x_S\in{]x_1^S,x_2^S[}\subset{\mathbb{R}}$ tel que $f(x_S)=y_S$.
\\
On admet $P(I)$ : $\exists{x^{(I)}}=x_S+x_1\cdot{o}+\dots{x_I\cdot{o^I}}$ tel que $\bar{f}(x^{(I)})=y^{(I)}$, avec
$y^{(I)}=y_S+y_1\cdot{o}+\dots{y_I\cdot{o^I}}$ et l'on d\'{e}montre $P(I+1)$.
\\
On cherche $x^{(I+1)}=x^{(I)}+x_{I+1}\cdot{o^{I+1}}$ tel que $\bar{f}(x^{(I+1}))=y^{(I+1)}$. On simplifie l'\'{e}quation et l'on trouve que $y_{I+1}=\bar{f}'(x^{(I)})\times{x_{I+1}}$. Le nombre $x_{I+1}$ est d\'{e}termin\'{e} si $\bar{f}'(x^{(I)})$ est inversible, c'est-\`{a}-dire $f'(x_S)\neq{0}$.  \qed

\paragraph{2} $\mathbb{R}_o$ n'est plus archim\'{e}dien mais on construit en 2. une extension intrins\`{e}que de $\mathbb{N}$ pour recouvrer cette propri\'{e}t\'{e}.
\\ \\
Le "co\^{u}t d'une r\'{e}volution non standard" est donc tr\`{e}s faible dans le cas de l'ANS* car, m\^{e}me s'il faut refaire toutes les d\'{e}monstrations classiques, elles sont toujours \'{e}l\'{e}mentaires et internes au cadre logique du "formalisme de Cauchy-Bolzano-Weierstrass-Cantor..." (E.Benoit) qui est le cadre g\'{e}n\'{e}ral des  \emph{El\'{e}ments de Math\'{e}matiques} de N.Bourbaki.
\\ \\
On pr\'{e}cise que l'\'{e}tude des propri\'{e}t\'{e}s diff\'{e}rentielles et int\'{e}grales des fonctions de $\mathbb{R}_o$ dans $\mathbb{R}_o$ ne s'appuie pas du tout sur le 4\`{e}me volume (F.V.R.) du Trait\'{e} de Bourbaki [11] puisque $\mathbb{R}_o$ n'est pas un corps valu\'{e} en tant qu'ensemble de scalaires, il n'est ni norm\'{e} ni complet en tant qu'espace vectoriel mais surtout la Topologie est la topologie induite par la relation d'ordre total et pas celle de l'Espace vectoriel norm\'{e}.
\par Pour mener \`{a} bien cette \'{e}tude, on \'{e}tudie d'abord dans le cadre g\'{e}n\'{e}ral de la Topologie bourbachique, les propri\'{e}t\'{e}s des prolongements analytiques des fonctions $C^{\infty}$ ou analytiques (cf. 1.2., 1.4 et 3.1) puis on d\'{e}finit des fonctions r\'{e}guli\`{e}res qui ne diff\`{e}rent des fonctions pr\'{e}c\'{e}dentes que sur un point de d\'{e}tail ($F(\mathbb{R})\subseteq{\mathbb{R}}$ ou non) et les preuves des Th\'{e}or\`{e}mes sont exactement les m\^{e}mes (cf. 1.3., 3.2 et 3.3).
\\ \\
Toute la section 5 appartient \'{e}galement au cadre de la Topologie bourbachique puisqu'il s'agit d'\'{e}tudier les propri\'{e}t\'{e}s alg\'{e}briques d'ensembles de nombres, totalement ordonn\'{e}s et de ce fait, topologiques (cf. aussi la \textbf{Remarque 5.6}).

\newpage

\section{De nouveaux nombres "entiers" d\'{e}finis \`{a} l'unit\'{e} pr\`{e}s.
\\Propri\'{e}t\'{e}s des ensembles $\mathbb{N}[\Sigma]$ et $\aleph$}

2.1. Conditions g\'{e}n\'{e}rales pour construire un ensemble de nombres entiers.
\\
2.2. Construction formelle de $\mathbb{N}[\Sigma]$ et $\mathbb{R}_o^{1+}$.
\\
2.3. Deux mod\`{e}les non standard de l'Arithm\'{e}tique de Peano.
\\
2.4. Le point de vue de I.Newton sur les nombres entiers infiniment grands.
\\ \\
\Large{\textbf{3. Une d\'{e}monstration du Th\'{e}or\`{e}me Fondamental de l'Analyse Non Standard}}
\\
\normalsize{3.1. R\'{e}solution de l'\'{e}quation du premier ordre pour le prolongement analytique d'une fonction analytique standard.
\\
3.2. G\'{e}n\'{e}ralisation aux fonctions r\'{e}guli\`{e}res, ou newtoniennes.
\\
3.3. R\'{e}solution d'une \'{e}quation diff\'{e}rentielle d'ordre quelconque pour des fonctions r\'{e}guli\`{e}res.}
\\ \\
\Large{\textbf{4. Le paradis newtonnien et la r\'{e}habilitation posthume de G.W.Leibniz}}
\\
\normalsize{4.1. Introduction \`{a} cette \'{e}tude d'Histonique.
\\
4.2. Deux caract\'{e}ristiques pr\'{e}alables de l'\'{e}difice newtonnien.
\\
4.3. Une modernisation du calcul des fluxions de Newton.
\\
4.4. Une comparaison des oeuvres de Newton et Leibniz.
\\
4.5. La Th\'{e}orie de l'Int\'{e}gration de I.Newton.
\\
4.6. Aux limites de l'\'{e}difice newtonien.
\\
4.7. Les faiblesses de l'\'{e}difice Leibnizien.
\\
4.8. La Th\'{e}orie de l'Int\'{e}gration de G.W.Leibniz.
\\
4.9. Conclusion de cette \'{e}tude d'Histonique.}
\newpage

\subsection{Conditions g\'{e}n\'{e}rales pour construire un ensemble de nombres entiers}
A la mani\`{e}re de Cantor [12, 13] et de Bourbaki [14]  on "construit" un nouvel ensemble de nombres entiers en quotientant un certain type d'ensembles par une certaine relation d'\'{e}quivalence (c'est l'\'{e}quipotence entre deux ensembles presque quelconques pour le "cardinal" et une bijection croissante entre deux ensembles bien ordonn\'{e}s pour l'"ordinal").
\\ \\
Deux autres conditions, \'{e}galement cantoriennes, semblent n\'{e}cessaires pour pouvoir consid\'{e}rer chaque classe d'\'{e}quivalence comme un entier naturel :
\\
1) L'Espace-quotient est bien et totalement ordonn\'{e}, il commence par $0$ suivi de tous les entiers standard.
\\
2) Il est muni d'une loi d'addition qui lui conf\`{e}re une structure de semi-groupe (mono\"{i}de) commutatif ou non.
\\ \\
La propri\'{e}t\'{e} de l'Espace-quotient d'\^{e}tre muni d'une application successeur qui lui conf\`{e}re la structure d'un "mod\`{e}le non standard de l'Arithm\'{e}tique de Peano" n'est pas \emph{n\'{e}cessaire} pour pouvoir consid\'{e}rer ses \'{e}l\'{e}ments comme des nombres entiers (elle n'est pas v\'{e}rifi\'{e}e par l'ensemble des nombres cardinaux ou ordinaux de Cantor).

\subsection{Construction formelle de $\mathbb{N}[\Sigma]$ et $\mathbb{R}_o^{1+}$}
Soit $\mathbb{R}_o^1=\lbrace{x_1=t+k\cdot{o}}/t\in{\mathbb{R}};k\in{\mathbb{Z}}\rbrace$.
\\
$J$ est l'ensemble des "intervalles" de $\mathbb{R}_o^1$, not\'{e}s $[[x_1,x'_1]]_1=$
\begin{center}
$[[x_1,x'_1+o[[_1=[x_1,x'_1]\cap{\mathbb{R}_o^1}$ si $x_1, x'_1\in{\mathbb{R}_o^1}$.
\end{center}

\begin{rem} Ni $\mathbb{R}_o$ ni $\mathbb{R}_o^1$ ne v\'{e}rifient la propri\'{e}t\'{e} de la borne sup\'{e}rieure dans leur propre topologie. Les intervalles de $J$ sont par contre tous born\'{e}s par leurs extr\'{e}mit\'{e}s dans la topologie de $\mathbb{R}_o^1$. Les bornes de ces intervalles sont donc uniques et bien d\'{e}finies.
\end{rem}
La relation d'\'{e}quivalence sur $J$ est la suivante :
\begin{center}
$[[x_1,x'_1]]_1\equiv{[[y_1,y'_1]]_1}$ ssi $x'_1-x_1=y'_1-y_1$.
\end{center}

On montre principalement que ces deux ensembles sont \textbf{deux mod\`{e}les non standard isomorphes de l'Arithm\'{e}tique de Peano}. Ce r\'{e}sultat va permettre d'intercaler entre $\mathbb{R}^+$ et
$\mathbb{R}_o^+$, un ensemble inductif $\mathbb{R}_o^{1+}$ ce qui permettra de faire des d\'{e}monstrations par induction dans $\mathbb{R}^+\subset{\mathbb{R}_o^{1+}}$.
\par
On pourra compter exactement le nombre de pas $o$ entre deux r\'{e}els standard positifs mais \textbf{on ne peut pas d\'{e}finir un successeur dans $\mathbb{R}^+$} car ces nombres de pas sont toujours infiniment grands comme on va le d\'{e}montrer. Ils n'ont donc pas de plus petite valeur possible.
\par
Tout cela sera plus clair apr\`{e}s l'introduction des notations $\otimes$ ("croix") et $\oslash$ ("slash").
\\ \\
$L$ est \underline{la} classe de l'intervalle $[[1\cdot{o}, 2\cdot{o}, 3\cdot{o},... x_1]]_1$.
\\
On l'\'{e}crit $L=\#[[o,x_1]]_1$ et l'on dit que $L$ est le \textbf{nombre d'\'{e}l\'{e}ments} de cette "suite" arithm\'{e}tique de premier terme et de raison $o$.
\\
On note $\Sigma=\#[[o,1]]_1$ et, puisque $k=\#[[o,k\cdot{o}]]_1$, par g\'{e}n\'{e}ralisation on note $x_1=L\otimes{o}$. Par cons\'{e}quent $k\cdot{o}=k\otimes{o}$ et
\begin{center}
$\Sigma\otimes{o}=1$.
\end{center}
$L$ ne d\'{e}pend que de $x_1$. On peut donc aussi noter $L=\Sigma\oslash{x_1}$ et puisque $1=\#[[o]]_1$, on a aussi :
\begin{center}
$\Sigma\oslash{o}=1$.
\end{center}

\begin{rem} Les deux \'{e}galit\'{e}s $\Sigma\otimes{o}=1$ et $\Sigma\oslash{o}=1$ ne signifient pas \emph{du tout} que $\Sigma$ et $o$ sont inverses l'un de l'autre puisque $\otimes$ et $\oslash$ ne sont pas des lois internes mais seulement des notations bien d\'{e}finies.
\end{rem}
Ces deux notations sont suffisantes pour d\'{e}montrer le Th\'{e}or\`{e}me suivant.

\begin{thm} Les ensembles $\mathbb{R}_o^{1+}$ et $\mathbb{N}[\Sigma]=J/\equiv$ sont en bijection par les applications bien d\'{e}finies $\varphi:\mathbb{R}_o^{1+}\longrightarrow{\mathbb{N}[\Sigma]}$, $x_1\longmapsto{\Sigma\oslash{x_1}}$ et
$\psi:\mathbb{N}[\Sigma]\longrightarrow{\mathbb{R}_o^{1+}}$, $L\longmapsto{L\otimes{o}}$.
\end{thm}
Preuve : $L=\#[[o,x_1]]_1$ s'\'{e}crit \`{a} la fois $L=\Sigma\oslash{x_1}$ et $x_1=L\otimes{o}$. En les combinant, on obtient les identit\'{e}s remarquables :
\begin{center}
$(\Sigma\oslash{x_1})\otimes{o}=x_1$ et $\Sigma\oslash{(L\otimes{o})}=L$
\end{center}
qui signifient exactement :
\begin{center}
$\psi\circ\varphi=Id_{\mathbb{R}_o^{1+}}$ et $\varphi\circ\psi=Id_{\mathbb{N}[\Sigma]}$.
\end{center} \qed

\begin{cor} On transf\`{e}re la structure de semi-groupe totalement ordonn\'{e} de $(\mathbb{R}_o^{1+},+,\leq)$ vers $(\mathbb{N}[\Sigma],\oplus,\preceq)$ par les applications $\varphi$ et $\psi$ et
\begin{center}
$L\oplus{M}=\Sigma\oslash{(L\otimes{o}+M\otimes{o})}$.
\end{center}
\end{cor}
On peut donc consid\'{e}rer (cf. 2.1) $\mathbb{N}[\Sigma]$ comme un nouvel ensemble de nombres entiers naturels.

\begin{prop} Les \'{e}l\'{e}ments de $\mathbb{N}[\Sigma]\setminus{\mathbb{N}}$ sont tous plus grands que tous les entiers standard de $\mathbb{N}$. Ils sont infiniment grands mais \emph{d\'{e}finis \`{a} l'unit\'{e} pr\`{e}s} (i.e. $L\oplus{1}\neq{L}$).
\end{prop}
Preuve : si $x_1$ n'est pas infinit\'{e}simal, $o\ll{x_1}$ et, par isomorphisme $L$ est sup\'{e}rieur \`{a} tous les entiers standard. Si $L\oplus{1}=L$ alors $\psi(L\oplus{1})=\psi(L)+o=\psi(L)$. Contradiction. \qed
\\ \\
Les entiers infiniment grands sont "d\'{e}finis" comme le sont les entiers finis standard. On peut d\'{e}finir maintenant une application "successeur".

\subsection{Deux mod\`{e}les non standard de l'Arithm\'{e}tique de Peano}
\begin{defn} Soient $s:\mathbb{R}_o^{1+}\longrightarrow{\mathbb{R}_o^{1+*}}$, $x_1\longmapsto{x_1+o}$ et
\\$S:\mathbb{N}[\Sigma]\longrightarrow{\mathbb{N}[\Sigma]^*}$, $L\longmapsto{L\oplus{1}}$. Ces applications "successeur" sont bien d\'{e}finies.
\end{defn}

\begin{thm} \textbf{(R\'{e}sultat principal)}
\\$(\mathbb{R}_1^+,s)$ et $(\mathbb{N}[\Sigma],S)$ sont deux mod\`{e}les non standard isomorphes de l'Arithm\'{e}tique de Peano [15, 16, 17].
\end{thm}
Preuve : l'application $s$ est bien une bijection. Il reste \`{a} prouver que $\mathbb{R}_o^{1+}$ est l'ensemble minimal qui contient $0$ et \textit{tous} ses successeurs par $s$. Soit $E\subseteq{\mathbb{R}_o^{1+}}$ tel que $0\in{E}$ et $s(E)\subseteq{E}$. On d\'{e}montre que $E=\mathbb{R}_o^{1+}$.
\par
Soit $x_1\in\mathbb{R}_o^{1+}$. $\Sigma\oslash{x_1}$ est le nombre d'\'{e}l\'{e}ments de $[[o,2\cdot{o},...x_1]]_1$ et l'on passe d'un terme \`{a} l'autre simplement en ajoutant $o$. Par cons\'{e}quent, $x_1$ est le $\Sigma\oslash{x_1}$-i\`{e}me successeur de $0$ par l'application $s$ et $x_1\in{E}$ puisque $0\in{E}$ et $s(E)\subseteq{E}$.
\par
Par isomorphisme, $(\mathbb{N}[\Sigma],S)$ est aussi un mod\`{e}le non standard de l'Arithm\'{e}tique de Peano. \qed

\begin{rem} Il fallait consid\'{e}rer $\Sigma\oslash{x_1}$ comme un nombre entier pour pouvoir appliquer la propri\'{e}t\'{e} $s(E)\subseteq{E}$ autant de fois (cf. 2.1).
\end{rem}
On peut d\'{e}finir une somme $\oplus$ ("plus") et un produit $\diamond$ ("fois") g\'{e}n\'{e}ralis\'{e}s dans le \textbf{prolongement inductif} $\aleph^+$ de $\mathbb{N}[\Sigma]^+$, par les \'{e}quations inductives [16, p.264] :
\begin{center}
$L\oplus{0}=L$ et $L\oplus{S(M)}=L\oplus{M}\oplus{1}$.
\\$L\diamond{1}=L$ et $L\diamond{S(M)}=L\diamond{M}\oplus{L}$.
\end{center}
Ces lois sont d\'{e}finies de proche en proche, par induction \hbox{dans $\aleph$.}
\\
On d\'{e}montre par induction toutes les propri\'{e}t\'{e}s du semi-anneau totalement ordonn\'{e} $(\aleph^+,\oplus,\diamond,\preceq)$ [14].

\begin{prop} On a $\aleph^+=$ \begin{center}
$\lbrace{\sum\limits_{0\leq{k}\leq{N}}a_k\Sigma^k}/(N\in{\mathbb{N}^{*}},a_0\in{\mathbb{Z}}
,$ $a_N>0)$ ou $(N=0, a_0\in\mathbb{N})\rbrace$
\end{center}
On note par un produit externe $a_k\Sigma^k$ le nombre entier $(\Sigma\oslash{a_k})\diamond{\Sigma}\diamond\dots\diamond{\Sigma}$ form\'{e} de $k-1$ produits.
\end{prop}
Preuve : $\aleph^+$ contient au moins tous ces nombres entiers puisque c'est un sur-anneau de $\mathbb{N}[\Sigma]=\varphi{(\mathbb{R}_o^{1+})}=$
\begin{center}
$\lbrace{\Sigma\oslash{a_1}\oplus{a_0}}/(a_1\in{\mathbb{R}^{+*}}$
 et $a_0\in{\mathbb{Z}})$ ou $(a_1=0$ et $a_0\in{\mathbb{N}})\rbrace$
\end{center}
 Il ne contient qu'eux puisque c'est l'ensemble minimal contenant $1$ et $\Sigma$, comme prolongement inductif de $\mathbb{N}[\Sigma]$. \qed

\subsection{Le point de vue de I.Newton sur les nombres entiers infiniment grands}
Sans souci de justification, I.Newton d\'{e}veloppe les m\^{e}mes id\'{e}es dans \underline{La m\'{e}thode des fluxions et des suites infinies} [6].
\par
Son traducteur en Fran\c{c}ais, le Marquis de Buffon exprime tr\`{e}s clairement cette conception tr\`{e}s originale pour nous aujourd'hui, de nombres infiniment grands mais d\'{e}finis \`{a} l'unit\'{e} pr\`{e}s, bien avant l'\'{e}poque cantorienne qui a fait pr\'{e}valoir l'id\'{e}e que les entiers infiniment grands ne sont pas d\'{e}finis \`{a} l'unit\'{e} pr\`{e}s.
\\ \\
Je cite longuement ce texte peu connu de 1740 :
\begin{quotation}
"\textit{Le Nombre n'est qu'un assemblage d'unit\'{e}s de m\^{e}me esp\`{e}ce ; l'unit\'{e} n'est point un nombre, l'unit\'{e} d\'{e}signe une seule chose en g\'{e}n\'{e}ral ; mais le premier Nombre 2 marque non seulement deux choses, mais encore deux choses semblables, deux choses de m\^{e}me esp\`{e}ce ; il en est de m\^{e}me de tous les autres Nombres.
\\ \\
Mais ces Nombres ne sont que des repr\'{e}sentations et n'existent jamais ind\'{e}pendamment des choses qu'ils repr\'{e}sentent ; les caract\`{e}res qui les d\'{e}signent ne leur donnent point de r\'{e}alit\'{e}, il leur faut un sujet, ou plut\^{o}t un assemblage de sujets \`{a} repr\'{e}senter pour que leur existence soit possible ; j'entends leur existence intelligible, car ils n'en peuvent avoir de r\'{e}elle ; or un assemblage d'unit\'{e}s ou de sujets ne peut jamais \^{e}tre que fini, c'est-\`{a}-dire, on pourra toujours assigner les parties dont il est compos\'{e}, par cons\'{e}quent le Nombre ne peut \^{e}tre Infini quelqu'augmentation qu'on lui donne.
\\ \\
Mais dira-t'on le dernier Terme de la suite naturelle 1, 2, 3, 4, etc n'est-il pas Infini ? n'y a-t-il pas des derniers Termes d'autres suites encore plus Infinis que le dernier terme de la suite naturelle ? Il para\^{\i}t que les Nombres doivent \`{a} la fin devenir Infinis, puisqu'ils sont toujours susceptibles d'augmentation ; \`{a} cela je r\'{e}ponds que \emph{cette augmentation dont ils sont susceptibles, prouve \'{e}videmment qu'ils ne peuvent \^{e}tre Infinis} ; je dis de plus que dans ces suites il n'y a pas de derniers Termes, que m\^{e}me leur supposer un dernier terme, c'est d\'{e}truire l'essence de la suite qui consiste dans la succession des Termes qui peuvent \^{e}tre suivis d'autres Termes et ces autres Termes encore d'autres, mais qui tous sont de m\^{e}me nature que les pr\'{e}c\'{e}dents, c'est-\`{a}-dire, tous finis, tous compos\'{e}s d'unit\'{e}s ; ainsi lorsqu'on suppose qu'une suite a un dernier Terme, et que ce dernier Terme est un nombre infini, on va contre la d\'{e}finition du nombre et contre la loi g\'{e}n\'{e}rale des suites"} [6, Pr\'{e}face, pages ix et x].
\end{quotation}

C'est exactement cette conception du nombre entier (fini ou non) de termes d'une suite arithm\'{e}tique de raison $o$, born\'{e}e des deux c\^{o}t\'{e}s, qui a \'{e}t\'{e} formalis\'{e}e en 2.2 et va \^{e}tre maintenant utilis\'{e}e pour une nouvelle Th\'{e}orie de l'Int\'{e}gration.

\newpage

\section{Une d\'{e}monstration du Th\'{e}or\`{e}me Fondamental de l'Analyse Non Standard}
On a montr\'{e} en 1.6 et en 2.4 que l'intuition math\'{e}matique de I.Newton exprim\'{e}e dans son \underline{Trait\'{e} des fluxions} [6] est rigoureusement correcte. On montre ici que l'intuition de G.W.Leibniz, portant sur une analogie entre Analyse Discr\`{e}te et Analyse Non Standard l'est \'{e}galement.

\begin{rem} 1) Dans l'article [5] on confirme le bien fond\'{e} de l'intuition physique de I.Newton concernant la gravitation, telle qu'elle est exprim\'{e}e dans les toutes derni\`{e}res lignes des \underline{Principia}.
\par 2) On appelle "Histonique" (cf. aussi 1.6) la tentative de reprendre des questions scientifiques l\`{a} o\`{u} les ont laiss\'{e}es les plus grands math\'{e}maticiens des si\`{e}cles pass\'{e}s, pour y donner de nouvelles r\'{e}ponses qui satisfont aux exigences de la Math\'{e}matique d'aujourd'hui (une \'{e}tude d'Histonique est donn\'{e}e \`{a} la section suivante).
\end{rem}

Conform\'{e}ment donc \`{a} l'intuition leibnizienne (cf. [10, p.208] et 1.6), tout se passe dans l'int\'{e}gration  d'une \'{e}quation diff\'{e}rentielle d'ordre quelconque, \textit{comme si} les sommes et les diff\'{e}rences \'{e}taient finies et d\'{e}finies, alors que les sommes sont infinies (elles ont un nombre infiniment grand mais d\'{e}fini de termes) et les diff\'{e}rences sont infinit\'{e}simales. \\La d\'{e}monstration de cette analogie intuitivement per\c{c}ue par G.W.Leibniz, entre Analyse Discr\`{e}te et Analyse Non Standard, devient possible parce que les ensembles $\mathbb{R}_o^{1+}$ et $\mathbb{R}_o$ ont "presque" les m\^{e}mes propri\'{e}t\'{e}s que les ensembles $\mathbb{N}$ et $\mathbb{R}$.

\subsection{R\'{e}solution de l'\'{e}quation du premier ordre pour le prolongement analytique d'une fonction analytique standard}
\begin{lem} (Th\'{e}or\`{e}me des sommes et diff\'{e}rences finies [10, p.204])
\\On d\'{e}finit sur l'ensemble $\mathcal{F}(\mathbb{N},\mathbb{R})$ (resp. $\mathcal{F}^*(\mathbb{N},\mathbb{R})$)
des suites num\'{e}riques standard (resp. de premier terme nul) deux op\'{e}rateurs $\Sigma$ et $\Delta$ qui sont \emph{r\'{e}ciproques l'un de l'autre} :
$$\Delta:\mathcal{F}^*(\mathbb{N},\mathbb{R})\rightarrow{\mathcal{F}(\mathbb{N},\mathbb{R})}, u\mapsto{v=\Delta{u}}$$
avec $v_n=u_{n+1}-u_n$.
$$\Sigma:\mathcal{F}^*(\mathbb{N},\mathbb{R})\rightarrow{\mathcal{F}(\mathbb{N},\mathbb{R})},
v\mapsto{u=\Sigma{v}}$$
avec $u_n=\sum\limits_{0\leq{k<n}}v_k$ ($u_0=0$).
\end{lem}
Preuve : on a $(\Delta\Sigma{v})_n=(\Sigma{v})_{n+1}-(\Sigma{v})_n=v_n$ et $(\Sigma\Delta{u})_n=\sum\limits_{0\leq{k}<n}(u_{k+1}-u_k)=u_n-u_0=u_n$
 apr\`{e}s $n$ simplifications. \qed
\\
On peut consid\`{e}rer chaque fonction $\bar{f}$ de $\mathbb{R}_o^{1}$ dans $\mathbb{R}_o$ comme form\'{e}e de deux "suites" $(\bar{f}(L\otimes{o}))_{L\in{\mathbb{N}[\Sigma]}}$ et $(\bar{f}(-L\otimes{o}))_{L\in{\mathbb{N}[\Sigma]}}$.

\begin{defn} Soit $f$ une fonction analytique en $t$, de $\mathbb{R}$ dans $\mathbb{R}$. Une primitive $G_1$ de $\bar{f}$ est une fonction de $\mathbb{R}_o^1$ dans $\mathbb{R}_o$ qui est une solution de l'\'{e}quation diff\'{e}rentielle
$$DG_1(x_1)=\bar{f}(x_1)\times{o}$$
pour tous les $x_1$ tels que $x_1^S\in{]t-R,t+R[}\subset{\mathbb{R}}$ o\`{u} $R$ est le rayon de convergence de la fonction analytique $f$ en $t$ ($R\in{\mathbb{R}^{+*}}$). De m\^{e}me pour $G$, fonction de $\mathbb{R}_o$ dans $\mathbb{R}_o$.
\end{defn}

On d\'{e}montre imm\'{e}diatement que si $f$ est analytique en $t$, de rayon de convergence $R\in{\mathbb{R}^{+*}}$, alors $\bar{f}$ est "analytique", c'est-\`{a}-dire $\bar{f}(t+x)=\sum\limits_{k\geq{0}}\frac{f^{(k)}(t)}{k!}x^k$ pour tous les $x\in{\mathbb{R}_o}$ tels que $|x_S|<R$.

\begin{thm} La primitive de $\bar{f}$ qui v\'{e}rifie la condition initiale $G_1(t)=a_0$ avec $a_0\in{\mathbb{R}_o}$, est la fonction de $\mathbb{R}_o^1$ dans $\mathbb{R}_o$ d\'{e}finie par la somme int\'{e}grale :
$$G_1(x_1)=a_0+\sum\limits_{y_L\in{[[t,x_1[[_1}}\bar{f}(y_L)\cdot{o}$$
si $x_1\geq{t}$ et
$$G_1(x_1)=a_0-\sum\limits_{y_L\in{[[-t+o,-x_1+o[[_1}}\bar{f}(-y_L)\cdot{o}$$
si $x_1\leq{t}$.
\end{thm}
Preuve : elle est la m\^{e}me que celle du Lemme, il y a seulement $\Sigma\oslash{x_1}$ simplifications. On peut aussi d\'{e}montrer la relation par induction dans $\mathbb{R}_o^{1+}$ et  $\mathbb{R}_o^{1-}$.
\\ \\
 \emph{Il faut maintenant pr\'{e}ciser les propri\'{e}t\'{e}s des fonctions int\'{e}grales $G$. On sait d\'{e}j\`{a} que ce ne sont pas des prolongements analytiques d'une fonction analytique standard car la condition $G_1(\mathbb{R})\subseteq{\mathbb{R}}$ n'a aucune raison d'\^{e}tre toujours v\'{e}rifi\'{e}e.}

\begin{thm} $G_1$, primitive de $\bar{f}$ valant $a_0$ en $t\in{\mathbb{R}}$, est une s\'{e}rie enti\`{e}re de $x_1\in{]t-R,t+R[}\subset{\mathbb{R}_o}$, i.e.
$$G_1(t+x_1)=a_0+\sum\limits_{l\geq{1}}\mathcal{A}_l\times{x_1^l}$$
avec $\mathcal{A}_l=\sum\limits_{m\geq{l-1}}\frac{f^{(m)}(t)}{m!}a_{m,l}\cdot{o^{m+1-l}}
\in{\mathbb{R}_o}$ pour certains coefficients standard $a_{m,l}$ que l'on pr\'{e}cisera.
\end{thm}
Preuve : 1) on d\'{e}montre la relation pour $x_1=k\cdot{o}$ avec $k\geq{0}$.
\\\hbox{$G_1(t+k\cdot{o})-G_1(t)=\sum\limits_{0\leq{n}<k}\bar{f}(t+n\cdot{o})\times{o}=$}
\\$\sum\limits_{0\leq{n}<k}\sum\limits_{m\geq{0}}\frac{f^{(m)}(t)}{m!}n^m\cdot{o^{m+1}}=
\sum\limits_{m\geq{0}}\frac{f^{(m)}(t)}{m!}\sum\limits_{0\leq{n}<k}p_m(n\cdot{o})\times{o}$.
\\ On note $p_m:\mathbb{R}_o\rightarrow{\mathbb{R}_0},x\mapsto{x^m}$ et $q_m$ la primitive de $p_m$ qui vaut $0$ en $0$. On admet provisoirement que $q_m$ est un polyn\^{o}me de degr\'{e} $m+1$ et l'on note $a_{m,l}$ le coefficient de $x^l\times{o^{m+1-l}}$ pour $1\leq{l}\leq{m+1}$ (cf. les deux Lemmes).
\\\hbox{$G_1(t+k\cdot{o})-G_1(t)=\sum\limits_{m\geq{0}}\frac{f^{(m)}(t)}{m!}\sum\limits_{0\leq{n}<k}
Dq_m(n\cdot{o})=\sum\limits_{m\geq{0}}\frac{f^{(m)}(t)}{m!}q_m(k\cdot{o})$}
$=\sum\limits_{m\geq{0}}\frac{f^{(m)}(t)}{m!}\sum\limits_{1\leq{l}\leq{m+1}}a_{m,l}
\cdot{(k\cdot{o})^l}\times{o^{m+1-l}}=
\sum\limits_{l\geq{1}}\mathcal{A}_l\times{(k\cdot{o})^l}$.
\\A peu pr\`{e}s le m\^{e}me calcul vaut pour $k\leq{0}$.

\par 2) le m\^{e}me calcul peut se faire pour $G_1(t+x_1)-G_1(t)$ avec $x_1=L\otimes{o}$ ou $-L\otimes{o}$ et $|x_1^S|<R$. On trouve de m\^{e}me que $G_1(t+x_1)=a_0+\sum\limits_{l\geq{1}}\mathcal{A}_l\times{x_1^l}$. \qed

\begin{lem}Soit $P_m:\mathbb{N}\rightarrow{\mathbb{N}}$ avec $P_m(N)=N^m$. On cherche $Q_m:\mathbb{N}\rightarrow{\mathbb{N}}$ telle que $\Delta{Q_m}(N)=Q_m(N+1)-Q_m(N)=P_m(N)$ et $Q_m(0)=0$. C'est la fonction $Q_m=\sum\limits_{1\leq{l}\leq{m+1}}a_{m,l}\cdot{P_l}$. Les coefficients $a_{m,l}$ sont tels que les matrices triangulaires d'ordre $m+1$ $(a_{m,l})_{1\leq{l}\leq{m+1}}$ et $(C_l^s)_{1\leq{l}\leq{m+1};0\leq{s}<l}$ soient inverses l'une de l'autre.
\end{lem}
Preuve : $Q_m=\sum\limits_{1\leq{l}\leq{m+1}}a_{m,l}P_l$ donne $\Delta{Q_m}=\sum\limits_{1\leq{l}\leq{m+1}}a_{m,l}\Delta{P_l}=P_m$ et
$\Delta{P_l}=\sum\limits_{0\leq{s}<l}C_l^sN^s$ donc
$\Delta{P_l}=\sum\limits_{0\leq{s}<l}C_l^sP_s$.
\qed

\begin{lem} Soit $p_m:\mathbb{R}_o^{1+}\rightarrow{\mathbb{R}_0}$, $p_m(x_1)=x_1^m$. On cherche la primitive $q_m$ de la fonction $p_m$ qui vaut $0$ en $0$. C'est le polyn\^{o}me $q_m=\sum\limits_{1\leq{l}\leq{m+1}}a_{m,l}[p_l\times{o^{m+1-l}}]$.
\end{lem}
Preuve : c'est la m\^{e}me que pr\'{e}c\'{e}demment. On montre que
$[p_m\cdot{o}]=\sum\limits_{1\leq{l}\leq{m+1}}a_{m,l}[Dp_l\cdot{o^{m+1-l}}]$ et
$[Dp_l\cdot{o^{m+1-l}}]=\sum\limits_{0\leq{s<l}}C_l^s[p_s\cdot{o^{m+1-s}}]$.\qed

\begin{rem} 1) On trouve par calcul que
$q_0=p_1$ ; $q_1=\frac{p_2}{2}-\frac{p_1}{2}o$ ; $q_2=\frac{p_3}{3}-\frac{p_2}{2}o+\frac{p_1}{6}o^2$ ; $q_3=\frac{p_4}{4}-\frac{p_3}{2}o+\frac{p_2}{4}o^2$ ; $q_4=\frac{p_5}{5}-\frac{p_4}{2}o+\frac{p_3}{3}o^2-\frac{p_1}{30}o^4$.
\par 2) On peut aussi exprimer les coefficients $a_{m,l}$ en fonction des Nombres de Bernouilli $B_p$. On trouve :
$$a_{m,l}=(-1)^{m+1-l}\sum\limits_{0\leq{p}\leq{m+1-l}}\frac{m!}{l!p!(m+1-l-p)!}B_p$$
\end{rem}

On prolonge la fonction $G_1$ sur $\mathbb{R}_o$ par $G(t+x)=G(t)+\sum\limits_{l\geq{1}}\mathcal{A}_l\times{x^l}$ pour tous les $x\in{\mathbb{R}_o}$ tels que $x_S\in{]t-R,t+R[\subset{\mathbb{R}}}$. On v\'{e}rifie que la fonction $G$ est une primitive de $\bar{f}$ sur $\mathbb{R}_o$ (la relation $DG=\bar{f}o$ est v\'{e}rifi\'{e}e sur $\mathbb{R}$ et l'on rappelle que la "d\'{e}rivation" se fait par rapport \`{a} la partie standard $t$ de $x$).

\begin{thm} $G$ est infiniment "d\'{e}rivable" (cf. 1.3) en $t$ et $\mathcal{A}_l=\frac{G^{(l)}(t)}{l!}$.
$$G(t+x)=\sum\limits_{l\geq{0}}\frac{G^{(l)}(t)}{l!}\times{x^l}$$
\end{thm}
Preuve : 1) On veut d\'{e}montrer qu'il existe $H$ tel que $(\forall{\varepsilon\in{\mathbb{R}_o^{+*}}})$
\\$(\exists{\eta\in{\mathbb{R}_o^{+*}}})
(\forall{x\in{\mathbb{R}_o}})|x|<\eta\Rightarrow{|G(t+x)-G(t)-H\times{x}|<\varepsilon\times{|x|}}$.
\\On prend bien s\^{u}r $H=\mathcal{A}_1$ et $b_2\in{\mathbb{R}^{+*}}$ tel que $|\mathcal{A}_2|<b_2$. \\L'ordre lexicographique donne $I=|G(t+x)-G(t)-\mathcal{A}_1\times{x}|=|\mathcal{A}_2\times{x^2}+\mathcal{A}_3\times{x^3}+
\dots|<b_2\cdot{x^2}<\varepsilon\times{|x|}$ d\`{e}s que $|x|<\eta<\frac{\varepsilon}{b_2}$. On prend $\varepsilon\ll{1}$ pour que $x$ soit infinit\'{e}simal.
\\ On a donc $G'(t)=\mathcal{A}_1$.

\par 2) On d\'{e}montre de m\^{e}me que \hbox{$G'(t+x)=G'(t)+\sum\limits_{l\geq{1}}(l+1)
\mathcal{A}_{l+1}\times{x^l}$}$=
G'(t)+\sum\limits_{l\geq{2}}l\mathcal{A}_l\times{x^{l-1}}$ pour $|x|<R$.
\\Il faut trouver $H$ tel que $I=|G(t+y)-G(t+x)-H\times{(y-x)}|<\varepsilon\times{|y-x|}$. C'est bien s\^{u}r $H=\sum\limits_{l\geq{1}}l\mathcal{A}_l\times{x^{l-1}}$.
\\$I=|\sum\limits_{l\geq{2}}\mathcal{A}_l[y^l-x^l-lx^{l-1}(y-x)]|=
|\mathcal{A}_2(y-x)^2+u^3|\leq{b_2\cdot{|y-x|^2}}$ si $\mathcal{A}_2^S\neq{0}$, $u^3$ d\'{e}signe un infinit\'{e}simal d'ordre 3 et $I<\varepsilon\times{|y-x|}$ dans les m\^{e}mes conditions que pr\'{e}c\'{e}demment.
\\Si $\mathcal{A}_2^S=0$, on remarque que $y^l-x^l-lx^{l-1}(y-x)=u_l|y-x|^2$ pour $p\geq{3}$ et un certain nombre infinit\'{e}simal $u_l$, alors $I\ll{|y-x|^2}<\varepsilon{|y-x|}$ pour $\eta=\varepsilon$.

\par 3) On pose $\mathcal{B}_l=(l+1)\mathcal{A}_{l+1}$ et \hbox{$G'(t+x)=G'(t)+\sum\limits_{l\geq{1}}\mathcal{B}_l\times{x^l}$} pour $|x|<R$ (cf. 2.).
\\Par r\'{e}currence, on d\'{e}montre que $G^{(k)}(t)=k!\mathcal{A}_k$ et que
\\\hbox{$G^{(k)}(t+x)=G^{(k)}(t)+
\sum\limits_{l\geq{1}}(l+1)(l+2)$}
$\dots{(l+k)}\mathcal{A}_{l+k}\times{x^l}$ (cf. 1. et 2.). \qed

Cette d\'{e}monstration n'utilise pas les valeurs particuli\`{e}res de $\mathcal{A}_l$ en fonction de $f$, elle est donc g\'{e}n\'{e}rale \`{a} toutes les s\'{e}ries enti\`{e}res de rayon de convergence standard non nul, fonctions qui vont \^{e}tre maintenant d\'{e}finies comme fonctions r\'{e}guli\`{e}res de $\mathbb{R}_o$ dans $\mathbb{R}_o$.

\begin{cor} $G(t+x+y)=\sum\limits_{k\geq{0}}\frac{G^{(k)}(t+x)}{k!}\times{y^k}$ pour \hbox{$|x_S|<R$}, $|y_S|<R$ et $|x_S+y_S|<R$. $G^{(k)}(t+x)=G^{(k)}(t)+\sum\limits_{l\geq{1}}\frac{G^{(l+k)}(t)}{(l+k)!}\times{x^l}$ pour $|x_S|<R$.
\end{cor}

\subsection{G\'{e}n\'{e}ralisation aux fonctions r\'{e}guli\`{e}res ou newtoniennes de $\mathbb{R}_o$ dans $\mathbb{R}_o$}
\emph{Les fonctions primitives $G$ n'ayant aucune raison de toujours v\'{e}rifier $G(\mathbb{R})\subseteq{\mathbb{R}}$, on d\'{e}finit des fonctions de $\mathbb{R}_o$ dans $\mathbb{R}_o$ qui ont les m\^{e}mes propri\'{e}t\'{e}s que les fonctions $\bar{f}$ (except\'{e}e $\bar{f}(\mathbb{R})\subseteq{\mathbb{R}}$) pour ensuite pouvoir int\'{e}grer des \'{e}quations diff\'{e}rentielles d'ordre quelconque.}

\begin{defn} Une fonction $F$ de $\mathbb{R}_o$ dans $\mathbb{R}_o$ est r\'{e}guli\`{e}re ou newtonienne en $t_o\in{\mathbb{R}}$ s'il existe une suite $(a_n)_{n\in{\mathbb{N}}}$ de nombres non standard telle que l'ensemble
\\$E=\{r\in{\mathbb{R}^{+*}}/(\forall{x\in{\mathbb{R}_o}})|x_S|<r\Rightarrow{F(t_o+x)=\sum\limits_{n\geq{0}}
a_n\cdot{x^n}\}}$
\\ne soit pas r\'{e}duit \`{a} l'ensemble $\emptyset$.
\end{defn}

\begin{rem} 1) $R$ est la borne sup\'{e}rieure de l'intervalle $E=]0,R[$ ou $]0,R]$ et $R>0$ avec la possibilit\'{e} que $R=+\infty$. $R$ est bien s\^{u}r le rayon de convergence de la "s\'{e}rie enti\`{e}re" en $t_0$.
\par 2) On ne peut recourir ici au Lemme d'Abel car la s\'{e}rie g\'{e}om\'{e}trique $\sum\limits_{n\geq{0}}x^n$ n'est pas "convergente" dans la topologie de $\mathbb{R}_o$ si $|x|<1$, seulement si $|x|\ll{1}$.
\end{rem}

\begin{prop} si $F$ est r\'{e}guli\`{e}re en $t_0$, $F$ est infiniment "diff\'{e}rentiable" en $t_0+x$ avec $|x_S|<R$. On note $F^{(n)}(t_0)=n!\cdot{a_n}$ et $F(t_0+x+y)=\sum\limits_{k\geq{0}}\frac{F^{(k)}(t_0+x)}{k!}\times{y^k}$ si $|x_S|<R$, $|y_S|<R$ et $|x_S+y_S|<R$ avec $F^{(k)}(t_0+x)=\sum\limits_{l\geq{0}}\frac{(k+l)!}{l!}a_{k+l}\times{x^l}$.
\end{prop}
Preuve : pour $x=0$, c'est la m\^{e}me que le \textbf{Th\'{e}or\`{e}me 3.9}.

\begin{rem} La derni\`{e}re formule n'est pas la d\'{e}riv\'{e}e d'une s\'{e}rie enti\`{e}re. On cherche comme aux points 2 et 3 de la preuve du Th\'{e}or\`{e}me pr\'{e}c\'{e}dent un nombre $H=a_k=\frac{F^{(k)}(t_0)}{k!}$ tel que $|F(t_0+x)-a_0-a_1\times{x}-\dots{-a_{k-1}\times{x^{k-1}}}-H\times{x^k}|<\varepsilon\times{|x|^k}$ pour $|x|<\eta$.
\end{rem}

\begin{defn} Une primitive d'une fonction quelconque $F$ est une fonction $G$ de $\mathbb{R}_o$ dans lui-m\^{e}me, qui est solution de l'\'{e}quation diff\'{e}rentielle $DG(t_0+x)=G(t_0+x+o)-G(t_0+x)=F(t_0+x)\times{o}$ pour $|x_S|<R$.
\end{defn}

\begin{rem} On montre en 4.6 que cette fonction est fortement ind\'{e}termin\'{e}e. C'est pourquoi on limite \`{a} partir de maintenant la notion Non Standard de "primitive" aux seules fonctions r\'{e}guli\`{e}res.
\end{rem}

\begin{thm} $F$ est une fonction r\'{e}guli\`{e}re de rayon de convergence $R\in{\mathbb{R}^{+*}}$ en $t_0\in{\mathbb{R}}$. Il existe une seule primitive de $F$ qui soit r\'{e}guli\`{e}re et vaille $0$ en $t_0$.
\\C'est la s\'{e}rie enti\`{e}re d\'{e}finie par $G(t_0+x)=\sum\limits_{l\geq{1}}\mathcal{A}_l\times{x^l}$ pour $|x|<R$  avec $\mathcal{A}_l=\sum\limits_{m\geq{l-1}}\frac{F^{(m)}(t_0)}{m!}a_{m,l}\cdot{o^{m+1-l}}=
\frac{G^{(l)}(t_0)}{l!}$.
\end{thm}
Preuve : 1) on montre comme pr\'{e}c\'{e}demment que la relation est v\'{e}rifi\'{e}e pour $x_1\in{\mathbb{R}_o^1}$ (cf. le \textbf{Th\'{e}or\`{e}me 3.5}).
\par 2) s'il existait deux prolongements r\'{e}guliers sur $\mathbb{R}_o$ qui co\"{\i}ncident sur $\mathbb{R}_o^1$, il existerait une fonction r\'{e}guli\`{e}re $K$ v\'{e}rifiant $K(t_0+x_1)=0$ pour tous les $x_1\in{\mathbb{R}_o^1}$ de $]t_0-R,t_0+R[$. On d\'{e}montre par r\'{e}currence que $K^{(p)}(t_0+x_1)=0$ en utilisant la condition de diff\'{e}rentiabilit\'{e} \`{a} l'ordre $p$ de $K$ ($p\geq{1}$) en $t_0+x_1$ dans $\mathbb{R}_o^1$ et $K(t_0+x)=K(t_0)+K'(t_0)\times{x}+\dots{=K(t_0)}=0$.
\par 3) la d\'{e}monstration de $G(t+x)=\sum\limits_{l\geq{1}}\frac{G^{(l)}(t)}{l!}\times{x^l}$ est la m\^{e}me que pr\'{e}c\'{e}demment (cf. le \textbf{Th\'{e}or\`{e}me 3.9}). \qed

\begin{rem} 1) On a : $$G'(t)=F(t)+\frac{1}{2}\cdot{F'(t)}\times{o}+\frac{1}{12}\cdot{F''(t)}\times{o^2}+
\dots{+\frac{a_{m,1}}{m!}}\cdot{F^{(m)}(t)}\times{o^m}
+\dots$$
Contrairement \`{a} l'Analyse Standard, la "d\'{e}riv\'{e}e" d'une fonction int\'{e}grale n'est pas exactement \'{e}gale \`{a} la fonction de d\'{e}part. C'est la "diff\'{e}rentielle" premi\`{e}re qui annule l'"int\'{e}gration" (cf. le \textbf{Th\'{e}or\`{e}me Fondamental}).
\par 2) D'apr\`{e}s la d\'{e}finition des fonctions r\'{e}guli\`{e}res, $F(t_0+u)=\sum\limits_{k\geq{0}}
\frac{F^{(k)}(t_0)}{k!}\times{u^k}$ pour $|u|\ll{1}$. La sous-section 1.5 se g\'{e}n\'{e}ralise donc sans changement aux fonctions r\'{e}guli\`{e}res de $\mathbb{R}_o$ dans $\mathbb{R}_o$.
\end{rem}

On \'{e}crit le Th\'{e}or\`{e}me Fondamental de l'Analyse Non Standard sous une forme analogue au \textbf{Lemme 3.2}.

\begin{thm} \textbf{Th\'{e}or\`{e}me Fondamental de l'Analyse Non Standard*}
\\$\mathfrak{F}(\mathbb{R}_o,\mathbb{R}_o)$ est l'ensemble des fonctions r\'{e}guli\`{e}res de $\mathbb{R}_o$ dans $\mathbb{R}_o$.
\\$\mathfrak{F}^*(\mathbb{R}_o,\mathbb{R}_o)$ est le sous-ensemble de $\mathfrak{F}(\mathbb{R}_o,\mathbb{R}_o)$ des fonctions valant $0$ en $0$.
\\On d\'{e}finit sur ces deux ensembles, deux op\'{e}rateurs $S$ et $D$ qui sont r\'{e}ciproques l'un de l'autre.
$$D:\mathfrak{F}^*(\mathbb{R}_o,\mathbb{R}_o)\rightarrow{\mathfrak{F}(\mathbb{R}_o,
\mathbb{R}_o)},G\mapsto{F}$$
$F$ est telle que $F(x)\times{o}=DG(x)=G(x+o)-G(x)$ et
$$S:\mathfrak{F}(\mathbb{R}_o,\mathbb{R}_o)\rightarrow{\mathfrak{F}^*(\mathbb{R}_o,
\mathbb{R}_o)},F\mapsto{G}$$
$G$ est l'unique primitive r\'{e}guli\`{e}re de $F$, valant $0$ en $0$ (cf. le Th\'{e}or\`{e}me pr\'{e}c\'{e}dent).
\end{thm}
Preuve : il suffit de v\'{e}rifier que $D$ est bien d\'{e}fini.
\\$DG=G'o+\frac{1}{2}G''o^2+\dots{=Fo}$, donc $F=G'+\frac{1}{2}G''o+\dots$.
\par 1) On montre que $F$ est "d\'{e}rivable" en $0$.
\\$I=|F(x)-F(0)-H\times{x}|=|[G'(x)+\frac{1}{2}G''(x)o+\dots{]}
-[G'(0)+\frac{1}{2}G''(0)o+\dots{]}-[G''(0)+\frac{1}{2}G'''(0)o+\dots{]}|=|\frac{1}{2}[G'''(0)+
\frac{1}{2}G^{iv}(0)o+\dots{]}+\frac{1}{6}[G^{iv}(0)+\frac{1}{2}G^v(0)o+\dots{]}x+
\dots{|}\times{x^2}$.
\\Si $G'''(0)_S\neq{0}$, on prend $\varepsilon$ infinit\'{e}simal, $b_2\in{\mathbb{R}^{+*}}$ tel que $|\frac{1}{2}G'''(0)_S|<b_2$ et $\eta{=\frac{\varepsilon}{b_2}}$. Si $G'''(0)_S=0$, $I=u\times{x^2}\ll{x^2}<\varepsilon{\times{x}}$ ($\eta=\varepsilon$).

\par 2) On a $F'(0)=H=G''(0)+\frac{1}{2}G'''(0)\times{o}+\dots$.
\\De m\^{e}me, on d\'{e}montre sur $\mathbb{R}_o$ que $F'=G''+\frac{1}{2}G'''o+\dots$ et, par r\'{e}currence sur $p$, que $F^{(p)}=G^{(p+1)}+\frac{1}{2}G^{(p+1)}o+\dots$.

\par 3) On d\'{e}montre que $F$ est r\'{e}guli\`{e}re par une permutation des signes de somme qui ne pose aucune difficult\'{e} \hbox{en ANS*.} $DG=F\times{o}$, donc $F(t+x)=\sum_{k\geq{0}}\frac{G^{(k+1)}(t+x)}{(k+1)!}\times{o^k}=
\sum_{p\geq{0}}\frac{1}{p!}[\sum_{k\geq{0}}
\frac{G^{(p+k+1)}(t)}{(k+1)!}o^k]\times{x^p}$ pour $|x_S|<R$ et $F(t+x)=\sum_{p\geq{0}}\frac{F^{(p)}(t)}{p!}\times{x^p}$

\par 4) On v\'{e}rifie que $DS=Id_{\mathfrak{F}(\mathbb{R}_o,\mathbb{R}_o)}$ et $SD=Id_{\mathfrak{F}^*(\mathbb{R}_o,\mathbb{R}_o)}$ (cf. le \textbf{Lemme 3.2}).
\qed

\begin{rem}On note $G=S[F\times{o}]$ lorsque \hbox{$G(x_1)=\sum\limits_{y_L\in{[[0,x_1[[_1}}F(y_L)\times{o}$}.
\end{rem}

\subsection{R\'{e}solution d'une \'{e}quation diff\'{e}rentielle d'ordre quelconque pour des fonctions r\'{e}guli\`{e}res}
On a r\'{e}solu l'\'{e}quation diff\'{e}rentielle du 1er ordre $DG=F\times{o}$ par $G=S[F\times{o}]$ (c'est le \textbf{Th\'{e}or\`{e}me 3.4}) dans le cas o\`{u} $F$ et $G$ sont r\'{e}guli\`{e}res. On fait de m\^{e}me pour int\'{e}grer une \'{e}quation diff\'{e}rentielle d'ordre sup\'{e}rieur.

\begin{thm} Soit une fonction $F$ quelconque de $\mathbb{R}_0$ dans $\mathbb{R}_0$. Toutes les fonctions $G_p$ qui sont solutions du syst\`{e}me diff\'{e}rentiel d'ordre p :
$D^kG_p(0)=C_k\times{o^k}$ pour $0\leq{k}<p$ et $C_k\in{\mathbb{R}_o}$ (ce sont les conditions initiales) et $D^pG_p(x_1)=F(x_1)\times{o^p}$, sont les fonctions $G_p$ d\'{e}finies sur $\mathbb{R}_o^1$ par
$$G_p(x_1)=S^p[F\times{o^p}]_{(x_1)}+C_0+C_1\times{x_1}+C_2\times{\frac{x_1(x_1-o)}{2}}
\dots{C_{p-1}}\times{B_{x_1}^{p-1}}$$
avec
$$
S^p[F\times{o^p}]_{(x_1)}=\sum\limits_{y_{L_p}\in{[[0,x_1[[_1}}\sum\limits_{y_{L_{p-1}}
\in{[[0,L_p\otimes{o}[[_1}}\dots
\sum\limits_{y_{L_1}\in{[[0,L_2\otimes{o}[[_1}}F(y_{L_p})\times{o^p}$$
et $B_{x_1}^k=\frac{x_1(x_1-o)\dots{(x_1-(k-1)o)}}{k!}$ (coefficient du bin\^{o}me dans $\mathbb{R}_o^{1+}$).
\end{thm}
Preuve : la d\'{e}monstration se fait par r\'{e}currence sur $p$. \\$G_p(x_1)-G_p(0)=\sum\limits_{y_{L_p}\in{[[0,x_1[[_1}}[DG_p(y_{L_p})-DG_p(0)]+
\sum\limits_{y_{L_p}\in{[[0,x_1[[_1}}DG_p(0)
\\=\sum\limits_{y_{L_p}\in{[[0,x_1[[_1}}
\sum\limits_{y_{L_{p-1}\otimes{o}}\in{[[0,y_{L_p}[[_1}}
D^2G_p(y_{L_2})+(\Sigma\oslash{x_1})\otimes{(C_1\times{o})}$.
\\Donc $S^2[D^2G_p]_{(x_1)}=S^2[D^2G_p-D^2G_p(0)]_{(x_1)}+
S^2[C_2\times{o^2}]_{(x_1)}=
S^3[D^3G_p]_{(x_1)}+\frac{\Sigma\oslash{x_1}(\Sigma\oslash{x_1}-1)}{2}\otimes{(C_2\times{o^2})}$ et
\\ $G_p(x_1)=S^3[D^3G_p]_{(x_1)}+C_0+C_1\times{x_1}+C_2\times{B_{x_1}^2}=
\dots{=S^p[F\times{o^p}]_{(x_1)}}+C_0+C_1\times{x_1}+\dots{C_{p-1}\times{B_{x_1}^{p-1}}}$.
\qed

\begin{defn} On appelle primitive $p$-i\`{e}me de la fonction $F$, toute solution de l'\'{e}quation diff\'{e}rentielle d'ordre $p$, $D^pG(x)=F(x)\times{o^p}$, o\`{u} $F$ est une fonction quelconque de $\mathbb{R}_0$ dans $\mathbb{R}_0$ (cf. \textbf{Remarque 3.16}).
\end{defn}

Il reste \`{a} donner la formule directe de $G_p(x_1)$ en fonction de $x_1$ pour montrer que son prolongement sur $\mathbb{R}_o$ est r\'{e}gulier \emph{lorsque $F$ est r\'{e}guli\`{e}re}.

\begin{lem} On sait calculer tous les coefficients $a_{m,l}^{(p)}$ de $q_m^{(p)}(x)=\sum\limits_{1\leq{l}\leq{m+p}}a_{m,l}^{(p)}\cdot{[x^l\times{o^{m+p-l}}]}$, primitive $p$-i\`{e}me de la fonction $p_m$ qui satisfait aux conditions initiales $D^kq_m^{(p)}(0)=0$ pour $0\leq{k}<p$.
\end{lem}
Preuve : 1) on trouve tout d'abord une relation portant sur les $a_{m,l}^{(p)}$ en Analyse Standard (cf. le \textbf{Lemme 3.6}).
\\On a $Q_m^{(p)}=\sum\limits_{1\leq{l}\leq{m+1}}a_{m,l}^{(p)}P_l$ avec $a_{m,m+p}^{(p)}=\frac{m!}{(m+p)!}$.
\\Donc $\Delta{Q_m^{(p)}}=\sum\limits_{1\leq{l}\leq{m+1}}a_{m,l}^{(p)}\sum\limits_{0\leq{k}\leq{l-1}}C_l^kP_k
=\sum\limits_{1\leq{l}\leq{m+1}}a_{m,l}^{(p)}\sum\limits_{1\leq{k}\leq{m+p}}C_l^{k-1}P_{k-1}$ avec $C_l^k=0$ d\`{e}s que $k\geq{l}$ (ici, $C_l^l=0$). On s'arrange pour que tous les indices varient entre 1 et $m+p$.
\\Alors $\Delta{Q_m^{(p)}}=\sum\limits_{1\leq{k}\leq{m+p}}b_{m,k-1}^{(p)}P_{k-1}$ avec
$[b_{m,k-1}^{(p)}]=(C_l^{k-1})[a_{m,l}^{(p)}]$, c'est le produit ligne par colonne d'une matrice triangulaire par une colonne.
\\De m\^{e}me, $\Delta^2{Q_m^{(p)}}=\sum\limits_{1\leq{k}\leq{m+p}}b_{m,n-1}^{(p)}\Delta{P_{k-1}}=
\sum\limits_{1\leq{n}\leq{m+p}}c_{m,n-1}^{(p)}P_{n-1}$ avec $[c_{m,n-1}^{(p)}]=(C_{k-1}^{n-1})[b_{m,k-1}^{(p)}]$.
\\Au final, on trouve que $\Delta^p{Q_m^{(p)}}=P_m=\sum\limits_{1\leq{n}\leq{m+p}}d_{m,n-1}^{(p)}P_{n-1}$ avec $[d_{m,n-1}^{(p)}]=(C_{k-1}^{n-1})^{p-1}(C_l^{k-1})[a_{m,l}^{(p)}]$.

\par 2) On d\'{e}montre de m\^{e}me pour $0\leq{s}\leq{m+p}$ que $\Delta^p{Q_s^{(p)}}=P_s=\sum\limits_{1\leq{n}\leq{m+p}}d_{s,n-1}^{(p)}P_{n-1}$ en posant $a_{s,l}=0$ d\`{e}s que $l>s+p$. On a toujours $[d_{s,n-1}^{(p)}]=(C_{k-1}^{n-1})^{p-1}(C_l^{k-1})[a_{s,l}^{(p)}]$.

\par 3) On \'{e}crit toutes ces relations sous la forme du produit matriciel $[P_{s-1}]=(d_{s,n-1}^{(p)})[P_{n-1}]$ et finalement : $$(a_{s,l}^{(p)})=[(C_{k-1}^{n-1})^{p-1}\times{(C_l^{k-1})}]^{-1}.$$ \qed

\begin{rem} On peut calculer par informatique, les premi\`{e}res valeurs de $a_{m,l}^{(p)}$.
\end{rem}

\begin{prop} L'unique fonction $G$ de $\mathbb{R}_o$ dans $\mathbb{R}_o$ qui v\'{e}rifie $D^kG(0)=0$ pour $0\leq{k}<p$ et $D^pG(x)=F(x)\times{o^p}$ est une s\'{e}rie enti\`{e}re, si $F$ est r\'{e}guli\`{e}re.
\end{prop}
Preuve : on sait que $G_p(x_1)=S^p[F\times{o^p}]_{(x_1)}$. $F$ est r\'{e}guli\`{e}re donc $G_p(x_1)=\sum\limits_{m\geq{0}}\frac{F^{(m)}(0)}{m!}\times{S^p[p_m\times{o^p}]_{(x_1)}}$ pour $|x_1^S|<R$ et $R>0$.
\\On conna\^{\i}t la primitive $p$-i\`{e}me de $p_m$, not\'{e}e $q_m^{(p)}$ qui v\'{e}rifie $D^kq_m^{(p)}(0)=0$ pour $0\leq{k}<p$ et \hbox{$S^p[p_m\times{o^p}]=S^pD^pq_m^{(p)}=q_m^{(p)}$.}
\\Donc $G_p(x_1)=\sum\limits_{l\geq{1}}\mathcal{A}_l^{(p)}\times{x_1^l}$ avec \hbox{$\mathcal{A}_l^{(p)}=\sum\limits_{m\geq{sup(l-p,0)}}\frac{F^{(m)}(0)}{m!}a_{m,l}^{(p)}
\cdot{o^{m+p-l}}$.}
De m\^{e}me pour $G$. \qed
\begin{prop} Cette s\'{e}rie enti\`{e}re est r\'{e}guli\`{e}re (cf. le \textbf{Th\'{e}or\`{e}me 3.9}).
\end{prop}

\newpage

{\scriptsize
\section{Le paradis newtonnien et la r\'{e}habilitation posthume de G.W.Leibniz}
\subsection{Introduction \`{a} cette \'{e}tude d'Histonique}
A propos de la querelle de priorit\'{e} qui a s\'{e}par\'{e} l'Angleterre du Vieux Continent sur un si\`{e}cle et demi, l'opinion de quelques Historiens des Sciences (le plus explicite est Ball [18, p.39-48]) est que G.W.Leibniz a eu connaissance d'une ou plusieurs lettres de I.Newton d\`{e}s 1675-1676, bien avant donc la publication de ses propres r\'{e}sultats, \`{a} partir de 1684. Il se serait donc "inspir\'{e}" des id\'{e}es de fluxions et de fluentes de Newton pour imaginer son propre symbolisme.
\\
On montre dans cette \'{e}tude d'Histonique jusqu'\`{a} quel point ces deux th\'{e}ories sont distinctes lorsqu'on les prolonge sous une forme compl\`{e}te et moderne : la fluxion $\dot{x}$ de Newton \emph{n'est pas} la d\'{e}riv\'{e}e $\frac{dy}{dx}$ de Leibniz, contrairement \`{a} ce qu'affirme un peu rapidement N.Bourbaki [10, p. 210]. Leibniz est ainsi r\'{e}habilit\'{e}, il n'a pas copi\'{e} sur Newton !
\\ \\
On ne rentrera pas ici dans le d\'{e}tail des reconstitutions \'{e}rudites de la pens\'{e}e de Newton [19, 20, 21], de Leibniz ou de Varignon [19] sauf quand elles servent d'arguments pour reconstituer ces deux \'{e}difices de Math\'{e}matiques. Il s'agit seulement ici de fournir une th\'{e}orie math\'{e}matique nouvelle, compl\`{e}te et simple, mais surtout "utile" pour les Math\'{e}matiques d'aujourd'hui, inspir\'{e}e principalement par l'oeuvre diff\'{e}rentielle et int\'{e}grale de I.Newton. D'une certaine mani\`{e}re c'est Newton qui a r\'{e}alis\'{e} le r\^{e}ve de Leibniz d'un Calcul Universel [10, p.210].
\par
On montre aussi que l'\'{e}difice newtonnien "engendre" math\'{e}matiquement celui de Leibniz mais le succ\`{e}s du dernier a fait oublier le premier, chronologiquement et \'{e}pist\'{e}mologiquement.

\subsection{Deux caract\'{e}ristiques pr\'{e}alables de l'\'{e}difice newtonnien}
Pour Newton, ni l'axe temporel ni les axes spatiaux ne sont gradu\'{e}s par des nombres (non standard), mais seulement par des points (cf. Fig. 1). La notation moderne $r_t$ ou $f(t)$ lui est donc \'{e}trang\`{e}re et c'est par une illusion r\'{e}trospective volontaire, que l'on va transcrire ses r\'{e}sultats de mani\`{e}re "moderne".
\\ \\
L'expos\'{e} moderne de cette th\'{e}orie math\'{e}matique nous oblige \`{a} renverser le point de vue qui a \'{e}t\'{e} celui de la sous-section 1.2 : au lieu de prolonger les fonctions $C^{\infty}$ standard sur les coupures infinit\'{e}simales, on part de l'id\'{e}e newtonienne que toutes les fonctions "physiques" (la distinction entre Math\'{e}matiques et Physique n'est pas encore faite \`{a} cette \'{e}poque, d'o\`{u} le titre de son Oeuvre ma\^{\i}tresse) sont NS*-continues \emph{\`{a} l'int\'{e}rieur des coupures infinit\'{e}simales}. On rappelle la d\'{e}finition de cette notion (cf. aussi 1.3).

\begin{defn} Une fonction de $\mathbb{R}_o$ dans $\mathbb{R}_o$ est NS*-continue [3] s'il n'est pas possible que
$$|x_2-x_1|\ll{|f(x_2)-f(x_1)|}.$$
\end{defn}

Ainsi, \`{a} l'int\'{e}rieur d'une coupure infinit\'{e}simale d'ordre $p$, la variation de $f$ est au plus infinit\'{e}simale d'ordre $q$, avec $q\geq{p}$. On a besoin de cette hypoth\`{e}se "naturelle" en Physique pour pouvoir d\'{e}finir les variations moyennes.

\subsection{Une modernisation du calcul des fluxions de Newton}

\paragraph{1} C'est dans les pages 190 \`{a} 193 des \textit{Principia} [22] que Newton est all\'{e} le plus loin dans l'usage de l'\'{e}criture alg\'{e}brique (cf. Fig. 1) que va privil\'{e}gier Leibniz (cf. 4.4.2). Partout ailleurs, il utilise les possibilit\'{e}s des dessins pour (se) repr\'{e}senter par des points les instants et par des intervalles les dur\'{e}es infinit\'{e}simales (cet effet de loupe des trac\'{e}s g\'{e}om\'{e}triques sera repris par les partisans de l'Analyse non standard de A.Robinson [2]).
\\
Apr\`{e}s avoir extrait "$DI=e-\frac{ao}{e}-\frac{nnoo}{2e^3}-\frac{anno^3}{2e^5}$ (aujourd'hui on \'{e}crirait $DI=e-\frac{a}{e}\cdot{o}-\frac{n^2}{2e^3}\cdot{o^2}-\frac{an^2}{2e^5}\cdot{o^3}$) de $DI^2=AQ^2-AD^2$ avec $AQ=n$ et $AD=a+o$, $n^2-a^2=e^2$, il nous dit (cf. Fig 1) :
\begin{quotation}
\textit{"Le premier terme qui est ici, $e$, repr\'{e}sentera toujours la longueur de l'ordonn\'{e}e $CH$ qui s'appuie sur le commencement de la quantit\'{e} ind\'{e}finie $o$... Le quatri\`{e}me terme d\'{e}termine la variation de la courbure ; le cinqui\`{e}me la variation de la variation, et ainsi de suite" }[22, p.192].
\end{quotation}

Ce qui s'\'{e}crit aujourd'hui : \\
Entre $C$ (d'abscisse $t\in{\mathbb{R}}$) et $D$ ($t+o$), la variation \emph{moyenne} de la fonction $r$ est $\dot{r}_t$ donn\'{e} par $Dr_t=r_{t+o}-r_t=\dot{r}_t\cdot{o}$ ($\dot{r}_t$ est en g\'{e}n\'{e}ral un nombre r\'{e}el non standard). Mais entre $D$ (d'abscisse $t+o$) et $E$ ($t+2o$), cette variation \emph{varie}.
\\
On consid\`{e}re s\'{e}par\'{e}ment la valeur pr\'{e}c\'{e}dente (qui se trouve de ce fait "conserv\'{e}e") de sa variation et
$$D\dot{r}_t=\dot{r}_{t+o}-\dot{r}_t=\ddot{r}_t\cdot{o},$$
$$D^2r_t=\dot{r}_{t+o}\cdot{o}-\dot{r}_t\cdot{o}=\ddot{r}_t\cdot{o^2}.$$

Donc, $\dot{r}$ (resp. $\ddot{r}$) est la variation moyenne d'ordre un (resp. deux) et non pas la d\'{e}riv\'{e}e premi\`{e}re (resp. seconde).
\\Newton ne va pas au del\`{a} de l'ordre 3. Sautons le pas et \'{e}crivons
$$D^kr_t=r^{[k]}_t\cdot{o^k}$$
avec $\dot{r}=r^{[1]}$ et $\ddot{r}=r^{[2]}$.

\paragraph{2} Newton savait (avec d'autres notations, cf. son Lemme 5 du Livre III) que $r_{t+k\cdot{o}}=r_t+C_k^1\cdot{Dr_t}+C_k^2\cdot{D^2r_t}+\dots{D^kr_t}$. Soit :
$$r_{t+u}=r_t+B_u^1\times{\dot{r}_t}+B_u^2\times{\ddot{r}_t}+\dots$$
pour n'importe quel $u=k\cdot{o}$ \'{e}l\'{e}ment de $[0[_1\subset{\mathbb{R}_o^{1+}}$.
\\
Voici la "s\'{e}rie formelle" la plus g\'{e}n\'{e}rale \`{a} laquelle I.Newton aurait pu parvenir.

\paragraph{3} Nous sommes tr\`{e}s loin de la Formule de Taylor (compatriote de Newton qui s'est en fait "inspir\'{e}" des r\'{e}sultats de Leibniz selon N.Bourbaki [10, p.209]).
\\
Une mani\`{e}re moderne de faire ce "passage" (qui n'est pas un "passage \`{a} la limite") de l'\'{e}difice \emph{discret} newtonnien \`{a} l'\'{e}difice \emph{continu} leibnizien, consiste \`{a} r\'{e}organiser les termes de la formule pr\'{e}c\'{e}dente en fonction des puissances de $u$ et l'on trouve
$$r_{t+u}=r_t+(\dot{r}_t+-\frac{1}{2}\ddot{r}_t\times{o}+\frac{1}{3}r^{[3]}_t\cdot{o^2}\dots{)}\times{u}$$
$$+\frac{1}{2}(\ddot{r}_t-r^{[3]}_t\times{o^2}
+\frac{11}{12}r^{[4]}_t\times{o^2}\dots{)}\times{u^2}+\frac{1}{6}(r^{[3]}_t
-\frac{3}{2}r^{[4]}_t\times{o}\dots{)}\times{u^3}+\dots.$$
Mais l'on ne sait plus, dans cette reconstitution math\'{e}matique de l'\'{e}difice newtonien que $r_t$ est, au sens moderne, ind\'{e}finiment d\'{e}rivable.

\paragraph{4} Comment alors sauter ce nouveau pas, qui s\'{e}pare Newton de Leibniz ?
\\
On note $a_n$ les termes entre les parenth\`{e}ses et $r_{t+u}=\sum\limits_{n\geq{0}}\frac{1}{n!}a_n\times{u^n}$ pour tous les \'{e}l\'{e}ments infinit\'{e}simaux $u$ de $\mathbb{R}_o^{1+}$. Cette s\'{e}rie formelle "converge" toujours dans $\mathbb{R}_o$, par d\'{e}finition. On appelle \emph{analytiques} les s\'{e}ries qui convergent en dehors de $[0[_1$ dans l'ensemble $\mathbb{R}_o^+$ et l'on \'{e}crit
$$r_{t+T}=\sum\limits_{n\geq{0}}\frac{1}{n!}a_n\times{T^n}$$
pour $|T|<R$ et $R\in{\overline{\mathbb{R}}^{+*}}$.
\\ \\
La d\'{e}monstration que cette fonction de la variable non standard $T$ est infiniment "diff\'{e}rentiable" et m\^{e}me r\'{e}guli\`{e}re, a \'{e}t\'{e} faite au \textbf{Th\'{e}or\`{e}me 3.9}. Alors $a_n=r^{(n)}(t)$ et $r^{(p)}(t+T)=\sum\limits_{n\geq{0}}\frac{1}{n!}a_{n+p}T^n$ (cf. le \textbf{Corollaire 3.10}). C'est le point maximal qu'aurait pu atteindre Leibniz s'il avait pu d\'{e}finir rigoureusement le symbole $\frac{dy}{dx}$  (cf. 4.7).

\paragraph{5} Le premier r\'{e}sultat de cette \'{e}tude d'Histonique est d'avoir clarifier la diff\'{e}rence d'approche entre ces deux initiateurs du Calcul Infinit\'{e}simal : approche \emph{discr\`{e}te} par pas de $o$ pour Newton, approche de la \emph{continuit\'{e}} et de la d\'{e}rivabilit\'{e} standard pour Leibniz ; variations moyennes $r^{[k]}$ d'un c\^{o}t\'{e}, variations instantan\'{e}es ou locales $r^{(k)}$ d'un autre.
\\ \\
On peut ainsi jeter aux oubliettes de l'Histoire la querelle de priorit\'{e} de cette immense d\'{e}couverte math\'{e}matique.

\paragraph{6} Les formules qui permettent de passer des valeurs moyennes aux d\'{e}riv\'{e}es sont les m\^{e}mes que celles des Propositions 1.29 et 1.31 :
\\
$F^{[p]}(x_1)=\sum\limits_{n\geq{p}}X^n_p\frac{F^{(n)}(x_1)}{n!}$ avec $X^n_p=\sum\limits_{k=0}^{p}(-1)^{p-k}C_p^kk^n$ et
\\
$\frac{F^{(n)}(x_1)}{n!}=\sum\limits_{p\geq{n}}(-1)^{p-n}K_{p-1}^{p-n}\frac{F^{[p]}(x_1)}{p!}o^{p-n}$ avec $A_k^p=\sum\limits_{n=1}^p(-1)^{p-n}K_{p-1}^{p-n}k^n$.

\subsection{Une comparaison des oeuvres de Newton et Leibniz}
\paragraph{1} Leibniz est r\'{e}habilit\'{e} mais son \'{e}difice ressemble beaucoup plus \`{a} l'immeuble moderne du Calcul des Variations en Analyse Standard. L'\'{e}difice newtonien est pour nous aujourd'hui beaucoup plus original et simple puisque la variation moyenne sur un intervalle d'amplitude infinit\'{e}simale permet d'\'{e}viter la question du "passage \`{a} la limite" qui ne sera vraiment r\'{e}solue qu'avec D'Alembert, au prix du formalisme des $\epsilon$ et des $\eta$.

\paragraph{2} Leibniz fait beaucoup moins de dessins que Newton dans ses articles sur le Nouveau Calcul et l'on pourrait presque s'en passer. Toute l'explication est port\'{e}e par l'\'{e}criture alg\'{e}brique et l'on peut dire que Leibniz est pour une \emph{Analyse analytique} alors que Newton est pour \emph{Analyse synth\'{e}tique} si l'on reprend les termes de l'\'{e}pisode Poncelet en G\'{e}om\'{e}trie [23, 24].
\\
Mais il utilise parfois des notations assez lourdes pour nous aujourd'hui car, pour lui, toute quantit\'{e} est positive alors que sa variation instantan\'{e}e peut \^{e}tre bien s\^{u}r positive ou n\'{e}gative. Il \'{e}crit donc $\pm{dy}$ ce que nous \'{e}cririons $dy$ et $\mp{dy}$ son oppos\'{e}.

\paragraph{3} En modernisant son propos, on comprend qu'il a laiss\'{e} \`{a} ses successeurs le soin de d\'{e}finir soigneusement  sa notation infinit\'{e}simale. Il se contente de dire qu'il cherche "la valeur de $dx:dy$ c'est-\`{a}-dire celle du rapport de $dx$ \`{a} $dy$" [19, p. 192].
\\
On \'{e}crira ici seulement $dy=y'_xdx$ car $dx$ n'est pas inversible dans $\mathbb{R}_o$. Leibniz a finalement b\'{e}n\'{e}fici\'{e} d'une facilit\'{e} de notation ($dx=\frac{dx}{dy}dy$) qui lui permet d'\'{e}crire $dx^{a}=a\cdot{x^{a-1}}dx$ ou $d\frac{1}{x^{a}}=-\frac{adx}{x^{a+1}}$ [19, p.180 ou 190].

\paragraph{4} En dehors de quelques formules tr\`{e}s g\'{e}n\'{e}rales et audacieuses  [10] comme "$d^2=dd$" (et m\^{e}me "$d^{-1}=\int$" qu'il \'{e}crit aussi "$d\int{x}$ aequ. $x$" et qui va devenir le Th\'{e}or\`{e}me Fondamental de l'Analyse), c'est P.Varignon essentiellement qui s'est le premier occup\'{e} du maniement des diff\'{e}rentielles secondes $d^2$, dans un contexte cin\'{e}matique.
\par
 M.Blay, sp\'{e}cialiste fran\c{c}ais de ce math\'{e}maticien fran\c{c}ais souligne en des termes tr\`{e}s mesur\'{e}s, une ambigu\"{\i}t\'{e} de son oeuvre qui reprend les termes de la distinction discret/continu mais en en inversant curieusement les destinataires (cf. 4.3.5) :
\begin{quotation}
\textit{"les difficult\'{e}s de la conceptualisation varignonienne r\'{e}sident principalement dans le fait que l'expression de l'accroissement de vitesse $dv$ et celle de la force $y$ impliquent si nous pouvons nous exprimer ainsi, une mod\'{e}lisation ambigu\"{e} du mode d'action de la force, en ce sens que celle-ci est cens\'{e}e agir... soit au tout premier instant... soit de fa\c{c}on constante et continue pendant tout l'intervalle de temps $dt$"} [19, p.206-7].
\end{quotation}

Cela peut s'exprimer math\'{e}matiquement ainsi :
\\
 Selon une conception \emph{discontinuiste} (\textit{"une force agit instantan\'{e}ment au d\'{e}but de cet intervalle de temps $dt$ puis n'agit plus jusqu'au d\'{e}but de l'intervalle de temps suivant" }[19, p.206]) de la force, la variation de vitesse est instantan\'{e}e : $v_{t_0^+}=v_{t_0}$ et $v_{t_0^-}=v_{t_0}+y\cdot{o}$. Donc $dv=y\cdot{dt}$ et la variation d'espace due \`{a} cette variation de vitesse est $d^2x=dv\cdot{dt}=y\cdot{dt^2}$.
\\
 Selon la conception \emph{continuiste} de Newton (\textit{"la force agirait effectivement continuellement et de fa\c{c}on constante pendant l'intervalle de temps $dt$, de telle sorte que l'accroissement de vitesse acquis \`{a} la fin de cet intervalle de temps soit encore $dv$" }[19, p.206]), \textit{"l'espace parcouru ne serait plus $ddx$ mais $\frac{1}{2}ddx$" }nous dit Michel Blay.

\begin{rem} Newton s'int\'{e}resse \`{a} la variation moyenne d'une fluente continue (et m\^{e}me infiniment d\'{e}rivable puisque c'est une s\'{e}rie enti\`{e}re convergente) par pas de temps discrets infinit\'{e}simaux (le plus souvent, $o$). Leibniz s'int\'{e}resse \`{a} la variation instantann\'{e}e d'une quantit\'{e} variable mais uniquement en ses valeurs standard.
\par
C'est pourquoi les deux peuvent \^{e}tre dit discontinuistes et continuistes selon que l'adjectif porte sur les propri\'{e}t\'{e}s des fonctions ou des ensembles de nombres ($\mathbb{R}$ et $\mathbb{R}_o^1$ sont discontinus dans $\mathbb{R}_o$ en un sens qui sera pr\'{e}cis\'{e} en 5.4).
\end{rem}

\paragraph{5} On voit \`{a} quel point les deux interpr\'{e}tations pr\'{e}c\'{e}dentes de l'acc\'{e}l\'{e}ration sont incompatibles mais seule l'interpr\'{e}tation continuiste  est "physiquement" juste car, pour Newton, la force ne peut sortir du n\'{e}ant. Il ne peut y avoir de discontinuit\'{e} dans les ph\'{e}nom\`{e}nes spatiaux-temporels, m\^{e}me dans les coupures infinit\'{e}simales secondes.
\par
Newton nous le dit, dans le Lemme X, Livre I, Section I des \textit{Principia} :
\begin{quotation}
\textit{"Les espaces qu'une force finie fait parcourir au corps qu'on presse, soit que cette force soit d\'{e}termin\'{e}e et immuable, soit qu'elle augmente ou diminue continuellement, sont dans le commencement du mouvement en raison doubl\'{e}e des temps" }[22].
\end{quotation}
M\^{e}me dans le cas simple o\`{u} l'acc\'{e}l\'{e}ration est "immuable", elle n'est pas \emph{\'{e}ternelle}, elle a un d\'{e}but et une fin et des solutions de raccordement par continuit\'{e} aux deux extr\'{e}mit\'{e}s. C'est l'hypoth\`{e}se de la N.S*-continuit\'{e}.

\paragraph{6} En r\'{e}sum\'{e}, la Th\'{e}orie de Leibniz para\^{\i}t beaucoup moins assur\'{e}e dans ses fondations que celle de Newton. C'est pourtant elle qui a gagn\'{e} au verdict de l'Histoire.

\subsection{La Th\'{e}orie de l'Int\'{e}gration de I.Newton}
\paragraph{1} C'est cette Th\'{e}orie qui a \'{e}t\'{e} d\'{e}velopp\'{e}e ici sous une forme moderne et compl\`{e}te mais en fait, Newton en est bien loin.
\\ \\
Dans le contexte cin\'{e}matique des \textit{Principia}, Newton montre qu'il est bien conscient de ce que le mouvement n'est que la somme infinie des d\'{e}placements infinit\'{e}simaux, ce qu'il exprime en ces termes :
\begin{quotation}
\textit{"Si d'un nombre \'{e}gal de particules on compose des temps quelconques \'{e}gaux [finis], les vitesses au commencement de ces temps seront comme les termes d'une progression continue pris \emph{par sauts}, en augmentant un nombre \'{e}gal de termes interm\'{e}diaires... Maintenant, soient diminu\'{e}es ces particules \'{e}gales de temps [jusqu'\`{a} devenir infiniment petites], et soit leur nombre augment\'{e} \`{a} l'infini, de sorte que l'impulsion de la r\'{e}sistance devienne \emph{continue} ; et les vitesses qui sont toujours en proportion continue dans les commencements des temps \'{e}gaux le seront encore dans ce cas"} [22, p.175].
\end{quotation}
On lira avec int\'{e}r\^{e}t la modernisation de la r\'{e}solution par Newton du probl\`{e}me de la r\'{e}sistance lorsqu'elle est proportionnelle \`{a} la vitesse dans [19, p.166-7].

\paragraph{2} Dans son \textit{Trait\'{e} des Fluxions} (cf. Fig. 2), il \'{e}tablit au Probl\`{e}me IX comment "Trouver l'Aire d'une Courbe propos\'{e}e quelconque" et pr\'{e}cise au Probl\`{e}me VII  \textit{"concevez que les Aires ACEB et ADB sont produites par le Mouvement des droites BD et BE le long de la ligne AB"} [6, p.86].
\\Si $z_t$ d\'{e}signe l'Aire ADB et $x_t=t$ celle du rectangle de c\^{o}t\'{e} $AC=1$, on a $\dot{x}_t=1$ et $Dz_t=\dot{z}_t\cdot{o}$ o\`{u} $\dot{z}_t$ repr\'{e}sente toujours la variation moyenne de l'aire entre $t$ et $t+o$. Alors, Newton dit simplement "par la relation donn\'{e}e des Fluxions, on trouve celle des Fluentes" [6, p.93].
\\ \\
A titre d'exemple un peu compliqu\'{e} (c'est l'Exemple 5), on montre comment Newton cherche $z$ en fonction de $x$ lorsque
$\dot{z}^3+a^2\dot{z}+ax\dot{z}-2a^3-x^3=0$.
\\
Newton sait que $z$ est une s\'{e}rie de puissances de $x$ (ou analytique), il cherche le terme constant $a$ puis remplace $\dot{z}$ par $a+px$, trouve $p=\frac{-1}{a}$, remplace $\dot{z}$ par $a-\frac{1}{a}x+px^2$, trouve $p=\frac{1}{64a}$, etc...

\paragraph{3} Newton est ainsi capable de trouver les premiers termes de toutes les fonctions $z$ telles que $\dot{z}$ soit une fonction alg\'{e}brique de $x$ mais il ne peut reconna\^{\i}tre les fonctions transcendantes Log, Arctg... Il ne peut que donner les valeurs approch\'{e}es des aires sans se pr\'{e}occuper si la s\'{e}rie formelle est ou non convergente en dehors de la coupure infinit\'{e}simale.
\\ \\
Par diff\'{e}rence, on a donn\'{e} en 3.1 les valeurs exactes (non standard) de toutes les sommes "int\'{e}grales" des valeurs prises par une fonction $\bar{f}\times{o}$ aux valeurs successives de $[[0,t[[_1\subset{\mathbb{R}_o^{1+}}$, m\^{e}me si on ne sait pas le plus souvent en calculer la partie standard.

\subsection{Aux limites de l'\'{e}difice newtonien}
On va r\'{e}pondre maintenant \`{a} une question plus "difficile" que s'est peut-\^{e}tre d\'{e}j\`{a} pos\'{e}e le lecteur :
\\
Pourquoi cette Th\'{e}orie Non Standard* de l'Int\'{e}gration a-t-elle \'{e}t\'{e} restreinte en 3.2 et 3.3 aux seules fonctions r\'{e}guli\`{e}res de $\mathbb{R}_o$ dans $\mathbb{R}_o$ (que l'on peut aussi nommer "newtoniennes") ?

\paragraph{1}
Il est possible de d\'{e}finir par induction dans $\mathbb{R}_o^{1+}$, une somme "int\'{e}grale" $G(x_1)=S[F\times{o}]_{(x_1)}$ pour une fonction $F$ \emph{quelconque} de $\mathbb{R}_o$ dans $\mathbb{R}_o$  mais il faut compter avec la forme ind\'{e}termin\'{e}e $\infty{\times}0$ comme va le montrer maintenant.
\\
Si $F=0$, toutes les fonctions $H$ qui sont constantes sur chaque coupure infinit\'{e}simale, v\'{e}rifient $DH(x_1)=0$ et
$$\sum\limits_{y_L\in{[[0,x_1[[_1}}0=\sum\limits_{y_L\in{[[0,x_1[[_1}}DH(y_L)=H(x_1)-H(0).$$

\paragraph{2} Cette somme  \emph{int\'{e}grale} $G$ est donc d\'{e}finie \`{a} une fonction $H$ pr\`{e}s telle que $H(x_1)=H(x_1^S)$, traduction du fait qu'on n'en conna\^{\i}t que les diff\'{e}rences premi\`{e}res.
\par
C'est encore pire avec les sommes multiples car $S^k[F\times{o^k}]_{(x_1)}=G(x_1)+a_1(x_1^S)x_1+a_2(x_1^S)x_1^2+\dots{+a_k(x_1^S)x_1^k}$, o\`{u} les fonctions de $x_1^S$ \`{a} valeurs dans $\mathbb{R}_o$ sont quelconques, traduction cette fois du fait qu'on ne conna\^{\i}t de cette fonction que ses diff\'{e}rences $k$-i\`{e}mes.

\paragraph{3} C'est pour \'{e}viter cette ind\'{e}termination de la somme int\'{e}grale que l'on s'est limit\'{e} aux seules fonctions newtoniennes. En effet, si $DH(x)=0$ partout, alors $H$ est une fonction constante (cf. le \textbf{Lemme}) et $S[0]_{(x_1)}=H(x_1)-H(0)=0$.

\begin{lem} Si $H$ est r\'{e}guli\`{e}re en $t_0$ et $DH=0$ sur $]t_0-R,t_0+R[$ alors $H$ est constante.
\end{lem}
Preuve : $H$ est "diff\'{e}rentiable" en $t_0+x$, $|x_S|<R$ donc $(\forall{\varepsilon\in{\mathbb{R}_o^{+*}}})(\exists{\eta\in{\mathbb{R}_o^{+*}}})
|DH(t_0+x)-H'(t_0+x)\times{o}|<\varepsilon\times{o}$ et $H'(t_0+x)=0$. De m\^{e}me, on d\'{e}montre par r\'{e}currence sur $p\geq{1}$ que $H^{(p)}(t_0+x)=0$ et $H(t_0+x)=H(t_0)+H'(t_0)\times{x}+\dots
{=H(t_0)}$.
\\ \\
La variation globale d'une fonction est alors la somme int\'{e}grale de ses variations infinit\'{e}simales (cf. la citation pr\'{e}c\'{e}dente de I.Newton).

\paragraph{4} Une autre justification de cette limitation volontaire de l'\'{e}difice newtonnien peut \^{e}tre donn\'{e}e : d'un point de vue physique, tous les mouvements spatiaux-temporels sont "newtonniens" (cf. 4.6.1 et 4.4.5).
\\
Il faut rappeler que, dans un contexte cin\'{e}matique, un objet quelconque ne peut changer de place de mani\`{e}re instantan\'{e}e, ni m\^{e}me sur une dur\'{e}e infinit\'{e}simale car toute vitesse moyenne doit rester finie.
\\
De m\^{e}me, si la vitesse semble changer de signe ou d'intensit\'{e} \`{a} l'instant d'un choc, des forces de contact dont l'intensit\'{e} varie de diverses mani\`{e}res en fonction de la vitesse, interviennent sur une dur\'{e}e infime mais finie (non infinit\'{e}simale).
\\
Enfin, une acc\'{e}l\'{e}ration constante appara\^{\i}t toujours de mani\`{e}re progressive.
\\ \\
Tous les "monstres" de l'Analyse Standard peuvent \^{e}tre ignor\'{e}s en Analyse Non Standard* car une fonction num\'{e}rique discontinue (prolong\'{e}e par des constantes dans les coupures infinit\'{e}simales) n'est pas NS*-continue. Toutes les fonctions non standard consid\'{e}r\'{e}es ici sont au moins infiniment "diff\'{e}rentiables" dans la topologie d'ordre de $\mathbb{R}_o$.

\subsection{Les faiblesses de l'\'{e}difice Leibnizien}
\paragraph{1} En se limitant au calcul des valeurs standard des d\'{e}riv\'{e}es, Leibniz a de fait \'{e}limin\'{e} du domaine des nombres, les \'{e}l\'{e}ments infinit\'{e}simaux. Ils ne sont plus que des auxilaires de calcul appel\'{e}s \`{a} dispara\^{\i}tre \`{a} la fin. C'est la m\^{e}me situation qu'ont connu un temps les nombres imaginaires [25].
\\ \\
Deux attitudes sont alors possibles :
\\
 Avec d'Alembert puis Cauchy, on \'{e}limine de mani\`{e}re habile de l'approche initiale de Leibniz ces \'{e}l\'{e}ments fant\^{o}mes. Toute l'Analyse Standard s'est construite sur ce refus des nombres infiniment petits et Leibniz en est d\'{e}j\`{a} en partie responsable.
\\
Avec Newton, on consid\`{e}re les nombres infinit\'{e}simaux comme autant r\'{e}els que les nombres standard et on apprend \`{a} calculer avec, ce que nous avons fait dans cet article (cf. aussi [5]).

\paragraph{2} Il y a un enjeu de taille derri\`{e}re ce choix d'accepter ou non les nombres r\'{e}els infiniment petits et les nombres entiers infiniment grands :
\\
Si la vitesse et l'acc\'{e}l\'{e}ration ne sont d\'{e}finies que pour les parties standard du temps, elles ne semblent plus agir que ponctuellement, par \`{a}-coups, \`{a} chaque coupure infinit\'{e}simale. Cette action \emph{discontinue} devient incompr\'{e}hensible du point de vue physique (cf. 4.4.4) mais elle est sans cons\'{e}quence puisque les \'{e}l\'{e}ments non standard de la droite num\'{e}rique ont disparu.
\par
Ainsi la Force, mais aussi tous les concepts que les Physiciens ont invent\'{e} depuis (le Travail, le Flux, etc...) sont d\'{e}finis sur la droite num\'{e}rique \emph{standard}, ils paraissent "continus" et m\^{e}me infiniment "d\'{e}rivables" en fonction du Temps, mais leurs modes d'action sont devenus incompr\'{e}hensibles.

\paragraph{3} Le lien \emph{r\'{e}el} des concepts physiques avec l'Espace-Temps-Masse a \'{e}t\'{e} irr\'{e}m\'{e}diablement perdu. Pour le r\'{e}tablir, un retour \`{a} Newton s'impose et une premi\`{e}re tentative de Physique Non Standard* est propos\'{e}e dans un article bient\^{o}t d\'{e}pos\'{e} sous ArchiV D.S. [5]. Il propose une premi\`{e}re compr\'{e}hension du "mode d'action" de l'interaction de gravitation \`{a} partir de sa premi\`{e}re formalisation math\'{e}matique, faite par I.Newton, il y a trois si\`{e}cles et demi.
\\ \\
 En faisant un clin d'oeil \`{a} K.Marx, on peut dire que "les Physiciens ont r\'{e}ussi \`{a} transformer le Monde, il s'agit maintenant de le comprendre".

\subsection{La Th\'{e}orie de l'Int\'{e}gration de G.W.Leibniz}
 Cette grande efficacit\'{e} des algorithmes de calcul de l'Analyse leibnizienne associ\'{e}e \`{a} une compl\`{e}te perte de "sens" se retrouve \`{a} l'identique dans la Th\'{e}orie de l'Int\'{e}gration de G.W.Leibniz. Un exemple d\'{e}taill\'{e} peut suffire \`{a} le montrer.

\paragraph{1} Afin de calculer l'aire du premier quadrant d'un disque unit\'{e}, Leibniz utilise la sous-tangente, c'est-\`{a}-dire le point d'intersection $T$ de la tangente en $M$ au cercle avec l'axe des $y$ (cf. Fig. 3).
\\
La consid\'{e}ration des triangles conduit \`{a} $\frac{y-z}{x}=\frac{1-x}{y}$, l'\'{e}quation du cercle \`{a} $dy=\frac{1-x}{y}dx$. L'aire $I$ est \'{e}gale \`{a} $1-\int\limits_{x=0}^{x=1}xdy$. Par changement de variable
$$1-I=\int\limits_{y=0}^{y=1}\frac{x(1-x)}{y}dx=\int\limits_{y=0}^{y=1}ydx
-\int\limits_{z=0}^{z=1}zdx=I-\int_{z=0}^{z=1}zdx.$$
Donc $I=\frac{1}{2}+\frac{1}{2}\int\limits_{z=0}^{z=1}zdx$. Une int\'{e}gration par parties donne
$I=1-\frac{1}{2}\int\limits_{x=0}^{x=1}xdz$.
\\
Les \'{e}quations pr\'{e}c\'{e}dentes donnent $x=\frac{2z^2}{1+z^2}$ que Leibniz d\'{e}veloppe en une s\'{e}rie de puissances ($x=2z^2-2z^4+2z^6\dots$) et int\`{e}gre terme \`{a} terme. Il obtient ainsi une valeur approch\'{e}e de $\frac{\pi}{4}$, ici $1-\frac{1}{3}+\frac{1}{5}-\frac{1}{7}\dots$.

\paragraph{2} On retrouve ici la puissance des algorithmes de calcul mais aussi la compl\`{e}te perte de sens. Comme pour les diff\'{e}rentielles et la diff\'{e}rence (cf. 4.4.3.), le lien entre int\'{e}grale et somme infinie de termes infinit\'{e}simaux est perdu.
\\
\textbf{L'id\'{e}e de diff\'{e}rence a disparu. L'id\'{e}e de somme a disparu aussi. Il ne reste plus chez Leibniz que la d\'{e}riv\'{e}e et la primitive (antid\'{e}riv\'{e}e) consid\'{e}r\'{e}es comme deux op\'{e}rateurs fonctionnels r\'{e}ciproques.}

\subsection{Conclusion de cette \'{e}tude d'Histonique}
Apr\`{e}s ce long d\'{e}tour par l'Histoire des Sciences, il semble \'{e}tabli que l'Analyse standard mais aussi la Physique Math\'{e}matique \'{e}l\'{e}mentaire souffrent de la perte d'une certaine "intuition" : les calculs se font bien s\^{u}r rigoureusement mais les liens entre diff\'{e}rentielle et diff\'{e}rence finie, int\'{e}grale et somme finie, entre concepts de Physique et Espace-Temps-Masse "r\'{e}el" [5] sont irr\'{e}m\'{e}diablement perdus.
\par
Bien s\^{u}r, un sens "analytique" est venu combler, avec le succ\`{e}s que l'on conna\^{\i}t, ce manque mais si l'on veut restaurer ces liens plus primitifs, le chemin \`{a} suivre est tout trac\'{e} : il faut consid\'{e}rer \`{a} nouveau les \'{e}l\'{e}ments de $\mathbb{R}_o$ comme aussi "r\'{e}els" que ceux de $\mathbb{R}$ et les \'{e}l\'{e}ments de $\mathbb{N}[\Sigma]$ ou $\aleph^+$ comme des entiers tout autant "naturels" que ceux de $\mathbb{N}$. Alors tout se passera comme si les pas de temps ($dt$, $o$, $u$) et les int\'{e}grales \'{e}taient des diff\'{e}rences et des sommes finies.
\par
Notre intuition \'{e}l\'{e}mentaire de l'Analyse discr\`{e}te (cf. le \textbf{Lemme 1.23} et le \textbf{Lemme 3.2}) pourra alors s'\'{e}tendre de l'infiniment petit \`{a} l'infiniment grand.
\\ \\
Ce projet devient r\'{e}alisable car $\mathbb{R}_o$ poss\`{e}de toutes les propri\'{e}t\'{e}s n\'{e}cessaires pour d\'{e}finir la diff\'{e}rentielle comme une diff\'{e}rence et $\aleph^+$ poss\`{e}de toutes les propri\'{e}t\'{e}s n\'{e}cessaires pour d\'{e}finir l'int\'{e}grale comme la somme d'un nombre d\'{e}fini de termes : $\mathbb{R}_o$ est une alg\`{e}bre (mais ce n'est pas un corps) topologique ; $\aleph^+$ est un mod\`{e}le non standard de l'arithm\'{e}tique de Peano et un semi-groupe additif.
\\ \\
On peut pr\'{e}ciser qu'il est possible de d\'{e}finir de bien d'autres mani\`{e}res une extension infinit\'{e}simale de $\mathbb{R}$ et une extension infinie-d\'{e}finie de $\mathbb{N}$.
\\Dans l'article [5], le symbole infinit\'{e}simal est not\'{e} $o_3$. C'est un \'{e}l\'{e}ment nilpotent (car $o_3^3=0$) de l'Alg\`{e}bre de Weil $\mathbb{R}_3=\mathbb{R}[Y]/(Y^3)$ [26]. Cette alg\`{e}bre totalement ordonn\'{e}e est tr\`{e}s diff\'{e}rente de celle \'{e}tudi\'{e}e ici puisque $\mathbb{R}_o$ ne poss\`{e}de pas d'\'{e}l\'{e}ments nilpotents.

\begin{rem} Seul $\mathbb{R}_o=\mathbb{R}[[X]]$ est intrins\`{e}que, c'est sa propri\'{e}t\'{e} sp\'{e}cifique importante. Il existe donc des nombres infiniment petits qui ne sont pas nilpotents [25].
\end{rem}

Si l'on veut comparer $o$ et $o_3$, il faut les plonger dans la sur-structure alg\'{e}brique et ordinale $\mathbb{R}[[X,Y]]/(Y^3)$. La seule fa\c{c}on de prolonger la relation d'ordre total sur $\mathbb{R}[[X]]$ et $\mathbb{R}[Y]/(Y^3)$ est de poser $0<o_3^2\ll{o_3}\ll{o}\ll{1}$ : les \'{e}l\'{e}ments nilpotents d'ordre 3 ($a\cdot{o_3}+b\cdot{o_3^2}$ avec $a\neq{0}$) ou 2 ($b\cdot{o_3^2}$ avec $b\neq{0}$) sont tous infiniment plus petits que tous les \'{e}l\'{e}ments infinit\'{e}simaux de $\mathbb{R}_o$.
\\ \\
Dans la correspondance de Leibniz, on trouve cette courte justification des nombres non standard :
\begin{quotation}
"Je suppose qu'il existe des quantit\'{e}s qui sont incomparablement plus grandes ou plus petites que d'autres"
\end{quotation}
Pour rendre r\'{e}el ce r\^{e}ve math\'{e}matique, il suffit de trouver une sur-structure alg\'{e}brique et ordinale de $\mathbb{R}_o$ et de $\aleph$. Cela sera fait au tout d\'{e}but de la section 5, avec le corps totalement ordonn\'{e} $(\Omega,+,\times,\leq)$.

\begin{rem} Dans le d\'{e}bat historique entre Nieuwentijt et Leibniz [26], il n'y a aucun vainqueur puisque, sans le savoir bien s\^{u}r, ils ne parlaient pas de la m\^{e}me structure : l'espace des nombres r\'{e}els de G.W.Leibniz pourrait s'apparenter \`{a} la structure de $\Omega$ alors que l'espace num\'{e}rique de Nieuwentijt pourrait \^{e}tre "simplement" $\mathbb{R}_2$.
\end{rem}
}
\newpage

\section{Principaux r\'{e}sultats sur les ensembles
\\de nombres}

5.1. Les extensions intrins\`{e}ques les plus simples de $\mathbb{N}$ et $\mathbb{R}$.
\\
5.2. Propri\'{e}t\'{e}s du corps totalement ordonn\'{e} $(\Omega,+,\times,\leq)$.
\\
5.3. Propri\'{e}t\'{e}s de l'alg\`{e}bre totalement ordonn\'{e}e \hbox{$(\Omega,+,\cdot,\times,\leq)$.}
\\
5.4. Propri\'{e}t\'{e}s ordinales de $(\Omega,\leq)$ et de $(\overline{\Omega},\leq)$.
\\
5.5. Une nouvelle caract\'{e}risation des structures $(\Omega,+,\times,\leq)$ et  $(\overline{\Omega},\leq)$.

\newpage
\subsection{Les extensions intrins\`{e}ques les plus simples de $\mathbb{N}$ et $\mathbb{R}$}
\begin{thm} $\mathbb{R}[[X]]=\lbrace{\sum\limits_{k\geq{0}}a_kX^k}/a_k\in{\mathbb{R}}\rbrace$ muni des lois usuelles et de l'ordre lexicographique tel que $0<X\ll{1}$, est une alg\`{e}bre int\`{e}gre totalement ordonn\'{e}e, not\'{e}e \hbox{$(\mathbb{R}_o,+,\cdot,\times,\leq)$}.
\end{thm}
Preuve : cf. la \textbf{Proposition 1.4}. \qed

\begin{prop} $(\mathbb{R}_o,+,\cdot)$ est un espace vectoriel topologique pour la topologie d'ordre.
\end{prop}
Preuve : on montre comme en Analyse standard, que la loi interne et la loi externe sont toutes les deux continues pour cette topologie. \qed

\begin{thm} $\aleph=\lbrace{\sum\limits_{0\leq{k}\leq{N}}a_kY^k}/(N\in{\mathbb{N}},a_0\in{\mathbb{Z},a_k\in{\mathbb{R}}})\rbrace$ muni de l'ordre lexicographique tel que $1\ll{Y}$ et des lois ad\'{e}quates, est un anneau commutatif unitaire et int\`{e}gre, totalement ordonn\'{e}, d'\emph{entiers tous d\'{e}finis \`{a} l'unit\'{e} pr\`{e}s}.
\end{thm}
Preuve : cf. la \textbf{Proposition 2.9}. On a ainsi facilement tous les nombres entiers n\'{e}gatifs, ils ont les m\^{e}mes propri\'{e}t\'{e}s que les nombres positifs. \qed

\begin{prop} $(\aleph^+,S)$ est un mod\`{e}le non standard de l'Arithm\'{e}tique de Peano qui prolonge d'une mani\`{e}re \emph{intrins\`{e}que} l'ensemble des entiers naturels $\mathbb{N}$.
\end{prop}
Preuve : c'est une formulation \'{e}quivalente au fait que $\aleph^+$ est le prolongement inductif des ensembles $\mathbb{N}$ et $\mathbb{N}[\Sigma]$ pour l'application $S$. Il ne d\'{e}pend d'aucun choix arbitraire. \qed

\begin{thm} $\Omega=\lbrace{\sum\limits_{k\leq{N}}a_kZ^k}/N\in{\mathbb{Z}},a_k\in{\mathbb{R},a_N\neq{0}}\rbrace\bigcup\{0\}$ muni des lois usuelles $+$, $\times$ et de l'ordre lexicographique tel que $1\ll{Z}$, est un \emph{corps} totalement ordonn\'{e}, archim\'{e}dien pour le prolongement intrins\`{e}que de $\mathbb{N}$.
\end{thm}
Preuve : $\aleph\subset\Omega$ pour $Y=Z=\Sigma$, $0\leq{k}$ et $a_0\in{\mathbb{Z}}$.
\\
$\mathbb{R}[[X]]\subset\Omega$ pour $X=Z^{-1}=o$, $N\leq{0}$.
\\
Un \'{e}l\'{e}ment non nul quelconque de $\Omega$ s'\'{e}crit $x=a_N\Sigma^N+a_{N-1}\Sigma^{N-1}+\dots{a_{-n}}o^n+\dots
=a_N\Sigma^N(1+\frac{a_{N-1}}{a_N}o+\dots{\frac{a_{-n}}{a_N}o^{N+n}}
+\dots)=a_N\Sigma^N(1+u)$, avec $u\in{\mathbb{R}_{o}}$ et $\vert{u}\vert\ll{1}$.
\par
Son inverse est $\frac{1}{a_N}o^N(1-u+u^2\dots{(-1)^nu^n}+\dots)\in\mathbb{R}_o$ si $N\geq{0}$ et $\frac{1}{a_N}\Sigma^{-N}(1-u+u^2\dots{(-1)^nu^n}+\dots)\in\Omega$ si $-N\geq{0}$ (cf. le \textbf{Corollaire 1.9}).
\\ \\
Soient $a\in\Omega^{*+}$ et $b\in\Omega$ avec $0<a<|b|$.
\par
On divise $b$ par $a$ et $x=\frac{b}{a}\in\Omega$. On appelle \emph{troncature enti\`{e}re} de $x=\sum_{k\leq{N}}c_kZ^k$, le nombre entier $L=[x]=\sum\limits_{k\leq{N}}d_kZ^k$ avec $d_k=0$ si $k<0$, $d_0=[c_0]$, partie enti\`{e}re de $c_0$ et $d_k=c_k$ si $k>0$.
\\
Par l'ordre lexicographique, $L\leq{x}<L+1$ et $La\leq{b}<(L+1)a$. Il y a toujours un multiple de $a$ qui d\'{e}passe le nombre $b$. \qed

\subsection{Propri\'{e}t\'{e}s du corps totalement ordonn\'{e} \hbox{$(\Omega,+,\times,\leq)$}}
\begin{rem} L'ensemble $\Omega$ n'est ni valu\'{e}, ni uniforme puisque
\begin{center}
$\vert{y-x}\vert\in{\Omega^+}$
\end{center}
or certains trait\'{e}s de N.Bourbaki [10, 27, 28] sont \'{e}crits pour un corps de scalaires $K$ \emph{norm\'{e}}, le plus souvent $\mathbb{R}$ ou $\mathbb{C}$.
\par
Le premier volume de Topologie [29] et la premi\`{e}re moiti\'{e} du second [30] visent par contre explicitement \`{a} "se d\'{e}barasser des nombres r\'{e}els" ([29], p.8) standard. Ils peuvent donc \^{e}tre ici utilis\'{e}s.
\end{rem}

\begin{prop} $(\mathbb{R}_o,+,\times)$ est un anneau topologique [30].
\end{prop}
Preuve : la preuve de la compatibilit\'{e} de l'addition avec la topologie d'ordre est classique. On sait que la compatibilit\'{e} de la multiplication est \'{e}quivalente aux deux axiomes $(AT_{IIIa})$ et $(AT_{IIIb})$ de [29, p.75].
$(AT_{IIIa})$ donne
$$(\forall{\varepsilon}\in\mathbb{R}_o^{+*})(\exists{\eta_1,\eta_2}\in\mathbb{R}_o^{+*})
(\forall{x,y}\in\mathbb{R}_o) \vert{x}\vert<\eta_1, \vert{y}\vert<\eta_2 \Longrightarrow{\vert{xy}\vert<\varepsilon}.$$
 On prend $0<o^{2n}<\varepsilon$ et $\eta_1=\eta_2=o^n$.
$(AT_{IIIb})$ donne
$$(\forall{\varepsilon}\in\mathbb{R}_o^{+*})(\exists{\eta}\in\mathbb{R}_o^{+*})
(\forall{x}\in\mathbb{R}_o) \vert{x}\vert<\eta \Longrightarrow{\vert{x_0x}\vert<\varepsilon}.$$
 Si $x_0^S\neq{0}$, on prend $\eta=\frac{\varepsilon}{\vert{x_0}\vert}$. Sinon, il suffit de prendre $\eta=\varepsilon$. \qed

\begin{prop}
$(\Omega,+,\times)$ est un corps topologique [30].
\end{prop}
Preuve : on sait d\'{e}j\`{a} que $\Omega$ est  un corps (cf. le \textbf{Th\'{e}or\`{e}me 5.5}). La continuit\'{e} \'{e}tant une propri\'{e}t\'{e} locale, la preuve que $(\Omega,+,\times)$ est un anneau topologique est la m\^{e}me que pr\'{e}c\'{e}demment. Il faut v\'{e}rifier l'axiome $(KT)$ de [30, p.83].
Soit, pour $x_0\neq{0}$,
$$(\forall{\varepsilon}\in\mathbb{R}_o^{+*})
(\exists{\eta}\in\mathbb{R}_o^{+*})
(\forall{x}\in\mathbb{R}_o) \vert{x_0}\vert<\eta \Longrightarrow{\vert{\frac{1}{x}-\frac{1}{x_0}}\vert<\varepsilon}.$$
 On prend $\eta=inf(\varepsilon\frac{x_0^2}{2},\frac{\vert{x_0}\vert}{2})$, alors $\vert{x}\vert>\frac{\vert{x_0}\vert}{2}$ et $\frac{\vert{x-x_0}\vert}{\vert{x}\vert\vert{x_0}\vert}<\frac{2\eta}{x_0^2}<\varepsilon$ m\^{e}me si $x_0=u$ est infinit\'{e}simal car $u\neq{0}$ est inversible. \qed

\begin{rem} N.Bourbaki reconna\^{i}t dans la Note historique de [30, p.223] que la droite num\'{e}rique standard est obtenue par lui par la compl\'{e}tion du groupe additif $\mathbb{Q}$ faite pour la premi\`{e}re fois par Cantor en 1872 (... et Meray en 1869 [7, p.245]) et non pas par la m\'{e}thode des coupures de Dedekind.
\\
La premi\`{e}re est classiquement faite dans un espace m\'{e}trique ou seulement uniforme (cf. la \textbf{Remarque 5.6}) mais une structure d'anneau totalement ordonn\'{e}e suffit.
\\
La seconde n'exige aucune structure alg\'{e}brique, une relation d'ordre total suffit pour d\'{e}finir la continuit\'{e} d'un ensemble (cf. 5.4.).
\end{rem}

\subsection{Propri\'{e}t\'{e}s de l'alg\`{e}bre totalement ordonn\'{e}e $(\Omega,+,\cdot,\times,\leq)$}
\begin{defn} Dans l'anneau $\mathbb{R}_o$ muni de sa topologie d'ordre, une suite d'\'{e}l\'{e}ments est une suite de Cauchy ssi
$$
(\forall{\varepsilon{\in\mathbb{R}_o^{+*}}})
(\exists{N}>0)(\forall{p}\geq{N})(\forall{q}\geq{N}) -\varepsilon<u_p-u_q<\varepsilon.
$$
\end{defn}

\begin{prop} $(\mathbb{R}_o,+,\times,\leq)$ est complet dans sa topologie d'ordre, i.e. toute suite de Cauchy d'\'{e}l\'{e}ments de $\mathbb{R}_o$ est convergente dans \hbox{$\Omega$ [31, p.5].}
\end{prop}
Preuve : soit une suite de Cauchy dans $\mathbb{R}_o$. Alors
$$(\forall{n>0})(\exists{N_n>0})(\forall{p>N_n}) u_{N_n}-o^n<u_p<u_{N_n}+o^n.$$
Ses $n$ premiers moments sont constants \`{a} partir du rang $N_n$. La suite $(u_{N_n})_{n\in\mathbb{N}}$ est donc convergente vers un nombre unique $l\in{\mathbb{R}_o}$ et
$$(\forall{\varepsilon}\in\mathbb{R}_o^{+*})(\exists{N_n}>0)(\forall{p>N_n}) \vert{u_p-l}\vert{<o^{n-1}}<\varepsilon.$$ \qed

\begin{prop} $(\Omega,+,\times,\leq)$ est complet dans sa topologie d'ordre.
\end{prop}
Preuve identique.

\begin{rem} $\mathbb{R}_o$ et $\Omega$ sont munis d'un "produit scalaire" qui v\'{e}rifie toutes les propri\'{e}t\'{e}s d'un produit scalaire standard, c'est le simple produit interne. On dira que ce sont des "Espaces de Hilbert" mais ce ne sont pas des Espaces de Hilbert.
\end{rem}

\subsection{Propri\'{e}t\'{e}s ordinales de $(\Omega,\leq)$ et de \hbox{$(\overline{\Omega},\leq)$}}
On montre que $\mathbb{R}$ n'est pas continu (au sens de Dedekind [31]) \emph{dans} l'ensemble $\mathbb{R}_o$. On cite tout d'abord R.Dedekind :
\begin{quotation}
Mais en quoi consiste exactement cette continuit\'{e}? Tout tient dans la r\'{e}ponse \`{a} cette question, et c'est par elle seule que l'on obtiendra \textbf{un fondement scientifique pour l'investigation de \emph{tous} les domaines continus}.
\\Je trouve l'essence de la continuit\'{e}... dans le principe suivant:
\\"Si tous les points de la droite se divisent en deux classes telles que tout point de la premi\`{e}re classe se situe \`{a} gauche de tout point de la deuxi\`{e}me, alors il existe un et un seul point qui produit cette r\'{e}partition de tous les points en deux classes, cette coupure de la droite en deux parties" [32, p.19-20].
\end{quotation}
Cela peut se formaliser de la mani\`{e}re suivante.

\begin{defn} $\mathbb{E}$ est un ensemble totalement ordonn\'{e}. On dit que $(C_g,C_d)$ est une coupure de $\mathbb{E}$ ssi
\begin{center}
$C_g\cup{C_d}=\mathbb{E}$ et $(\forall{x}\in{C_g})(\forall{x'}\in{C_d}) x\leq{x'}$.
\end{center}
\end{defn}

\begin{rem} On voit plus loin pourquoi on permet aux deux parties de la coupure d'avoir un \'{e}l\'{e}ment en commun. Mais ce n'est pas une obligation (sinon, la propri\'{e}t\'{e} serait triviale).
\end{rem}

\begin{defn} $\mathbb{E}\subseteq\mathbb{F}$ sont deux ensembles totalement ordonn\'{e}s. $\mathbb{E}$ \emph{est continu dans} $\mathbb{F}$ ssi, quelle que soit la coupure $(C_g,C_d)$ de $\mathbb{E}$, il existe toujours \emph{un seul} $z\in\mathbb{F}$ tel que
\begin{center}
$(\forall{x}\in{C_g})(\forall{x'}\in{C_d})$ $x\leq{z}\leq{x'}$.
\end{center}
$\mathbb{E}$ est continu (discontinu) ssi il est (n'est pas) continu dans lui-m\^{e}me.
\end{defn}

\begin{rem} D'apr\`{e}s la d\'{e}finition pr\'{e}c\'{e}dente, R.Dedekind d\'{e}montre dans [32] que $\mathbb{Q}$ est discontinu mais continu dans $\mathbb{R}$ et que $\mathbb{R}$ est continu.
\end{rem}

\begin{prop} $\mathbb{R}$ n'est pas continu \emph{dans} l'ensemble $\mathbb{R}_o$. Il est discret dans la topologie d'ordre de $\mathbb{R}_o$.
\end{prop}
Preuve : toute une coupure infinit\'{e}simale vient s'intercaler entre les deux parties d'une coupure de Dedekind de la droite num\'{e}rique $\mathbb{R}$ standard et
$(\forall{t\in\mathbb{R}})]t-o;t+o[\cap\mathbb{R}=\{t\}$. \qed
\\ \\On met en relation la continuit\'{e} avec deux autres propri\'{e}t\'{e}s bien connues, la densit\'{e} et le Th\'{e}or\`{e}me de Bolzano.

\begin{defn} $\mathbb{E}\subseteq\mathbb{F}$ sont deux ensembles totalement ordonn\'{e}s. $\mathbb{E}$ \emph{est dense dans} $\mathbb{F}$ ssi
$$(\forall{x,y}\in{\mathbb{F}})(\exists{z}\in{\mathbb{E}})x<y\Rightarrow{x<z<y}.$$
$\mathbb{E}$ est dense ssi il est dense dans lui-m\^{e}me.
\end{defn}

\begin{lem} $(\mathbb{E},\leq)$ totalement ordonn\'{e}. $\mathbb{E}$ est continu ssi il v\'{e}rifie le Th\'{e}or\`{e}me de Bolzano et il est dense.
\end{lem}
Preuve : 1) Soit $E$ un ensemble dense v\'{e}rifiant, si $X\subset\mathbb{E}$ est major\'{e} (resp. minor\'{e}), alors $X$ est born\'{e} sup\'{e}rieurement (resp. inf\'{e}rieurement).
\\
Soit $(C_g,C_d)$ une coupure de $\mathbb{E}$, $C_g$ (resp. $C_d$) est major\'{e} (resp. minor\'{e}) donc born\'{e} par $b_g$ (resp. $b_d$). On a $b_g\leq{b_d}$. Si $b_g<b_d$, la condition $C_g\cup{C_d}=\mathbb{E}$ n'est plus satisfaite du fait de la densit\'{e} de $\mathbb{E}$, donc $b=b_g=b_d$. Par cons\'{e}quent, $$(\forall{x\in{C_g}})(\forall{y\in{C_d}}) x\leq{b}\leq{y}.$$
\par 2) R\'{e}ciproquement, si $(\mathbb{E},\leq)$ est continu et si $X\subset\mathbb{E}$ est major\'{e}. On montre que $C_g=\{x_g\in\mathbb{E}/(\exists{x}\in{X})x_g\leq{x}\}$ et $C_d=\{x_d\in\mathbb{E}/(\forall{x}\in{X})$ $x<x_d\}$ forment une coupure de $\mathbb{E}$ :
\\
$C_g\cup{C_d}=\mathbb{E}$ ; soit $x_g\in{C_g}$, $(\exists{x}\in{X})$ $x_g\leq{x}$ et $(\forall{x_d}\in{C_d})$ $x<x_d$ (la deuxi\`{e}me condition est satisfaite). Il existe donc $x_0$ unique plus petit des majorants de $X$ et $X$ est born\'{e} sup\'{e}rieurement.
\par
$\mathbb{E}$ est dense car si ce n'\'{e}tait pas le cas, il existerait $x_0<y_0$ tels que $]x_0,y_0[\cap{\mathbb{E}}=\emptyset$. Alors, $C_g=\{x\in{\mathbb{E}}/x\leq{x_0}\}$ et $C_d=\{y\in{\mathbb{E}}/y\geq{y_0}\}$ formeraient une coupure de $\mathbb{E}$ et il existerait deux nombres $z=x_0$ et $z=y_0$ tels que
$(\forall{x}\in{C_g})(\forall{y}\in{C_d})$ $x\leq{z}\leq{y}$. Contradiction.
\qed

\begin{prop} $(\mathbb{R}_o,\leq)$ n'est pas continu.
\end{prop}
Preuve : on sait d\'{e}j\`{a} que $\mathbb{R}_o$ ne v\'{e}rifie pas la propri\'{e}t\'{e} de la borne sup\'{e}rieure (cf. la \textbf{Remarque 2.1}). On peut aussi prendre $C_g=\mathbb{R}_o^-\cup{[0[}$ et $C_d=\mathbb{R}_o^+\backslash{[0[}$. On v\'{e}rifie les deux propri\'{e}t\'{e}s d'une coupure mais il n'existe pas $\epsilon\in{\mathbb{R}_o}$ tel que $\epsilon$ soit la borne sup\'{e}rieure de $[0[$. \qed
\\
\begin{quote}
    \textbf{Pour que la droite num\'{e}rique non standard v\'{e}rifie le Th\'{e}or\`{e}me de Bolzano, il suffit de consid\'{e}rer la droite num\'{e}rique achev\'{e}e $\overline{\mathbb{R}}=\mathbb{R}\cup{\{-\infty,+\infty\}}$ qui est un espace topologique totalement ordonn\'{e}} [30, Chap.IV, p.145 \`{a} 150, p.162] \textbf{m\'{e}trisable} [31, p.51] \textbf{.}
\end{quote}
On ajoute \`{a} l'ensemble $\mathbb{R}_o$ tous les points isol\'{e}s qui ferment les coupures infinit\'{e}simales.

\begin{defn} On pose $\overline{\mathbb{R}_o}=\mathbb{R}_o\cup{F}$ avec
$$F=\{t+\sum\limits_{1\leq{k}\leq{P}}a_k\cdot{o^k}/t\in\mathbb{R}, P\in\mathbb{N}^*, a_k\in\mathbb{R}, 1\leq{k}<P, a_P=\pm\infty\}.$$
\end{defn}

\begin{rem} $\overline{\mathbb{R}_o}$ est totalement ordonn\'{e} par l'ordre lexicographique, l'ordre sur $\overline{\mathbb{R}}$ et $\infty\cdot{o}\ll{1}$, i.e. $(\forall{n}\in{\mathbb{N}})$ $+\infty\cdot{o}<\frac{1}{n}$. On note $\epsilon=+\infty\cdot{o}$ et $\epsilon\cdot{o}=+\infty\cdot{o^2}$ etc... les bornes sup\'{e}rieures des coupures successives en 0 et l'on peut rigoureusement \'{e}crire dans $\overline{\mathbb{R}_o}$ que $]t[=[t-\epsilon,t+\epsilon]$.
\\ $\overline{\mathbb{R}_o}$ a perdu toutes ses propri\'{e}t\'{e}s alg\'{e}briques puisque l'addition n'est plus interne.
\end{rem}
La d\'{e}monstration compl\`{e}te du th\'{e}or\`{e}me suivant est un peu longue mais cela prouve que l'ajout des \'{e}l\'{e}ments infinis de $\overline{\mathbb{R}}$ \`{a} $\mathbb{R}_o$ lui comble toutes ses discontinuit\'{e}s.

\begin{thm} $(\overline{\mathbb{R}_o},\leq)$ est continu.
\end{thm}
{\scriptsize
Preuve : $\overline{\mathbb{R}_o}$ est dense dans lui-m\^{e}me.
\\Soit $X\subset{\overline{\mathbb{R}_o}}$ major\'{e}. On d\'{e}montre qu'il existe $b=\sum\limits_{i\geq{0}}b_io^{i}$ borne sup\'{e}rieure dans la topologie d'ordre de $\overline{\mathbb{R}_o}$. On appelle Troncature d'ordre I de $x=\sum\limits_{i\geq{0}}a_io^{i}$, le nombre $T_I(x)=\sum\limits_{0\leq{i}\leq{I}}a_io^{i}$ pour $I\in{\mathbb{N}}$. On d\'{e}montre par r\'{e}currence la proposition $P(I)$ : \\$T_I(\overline{\mathbb{R}_o})$ v\'{e}rifie la propri\'{e}t\'{e} de Bolzano, i.e. si $X\subset{\overline{\mathbb{R}_o}}$ est major\'{e}, alors $T_I(X)$ est born\'{e} sup\'{e}rieurement dans la topologie d'ordre de $T_I(\overline{\mathbb{R}_o})$.

\paragraph{1}
$P(0)$ : $T_0(\overline{\mathbb{R}_o})=\mathbb{R}$ v\'{e}rifie la propri\'{e}t\'{e} de la borne sup\'{e}rieure. On d\'{e}montre seulement $P(1)$, la g\'{e}n\'{e}ralisation est imm\'{e}diate.
\\
Soit $X\subset{T_1(\overline{\mathbb{R}_o})}$ major\'{e}, donc $T_0(X)$ est major\'{e} aussi et il est born\'{e} par $b_0=SupT_0(X)\in{\mathbb{R}}$.
\\
Soit $F:X\cap{]b_0[}\rightarrow\overline{\mathbb{R}}$, $x=b_0+a_1\cdot{o}\mapsto{a_1}$.
\\
$ImF\subset{\overline{\mathbb{R}}}$ donc il est born\'{e} sup\'{e}rieurement [30, p.34] par $b_1\in{\overline{\mathbb{R}}}$. On d\'{e}montre en m\^{e}me temps dans les trois cas ($b_1=-\infty,+\infty$ et $b_1$ fini) que $SupX=b_0+b_1\cdot{o}$ dans $T_1(\overline{\mathbb{R}_o})$.

\paragraph{1.1}
Soit $x\in{X}$, $x=a_0+a_1\cdot{o}$, $a_0\leq{b_0}$ ou $a_0=b_0$ et $a_1\leq{b_1}$ donc $x\leq{b_0+b_1\cdot{o}}$. On montre que c'est le plus petit majorant de $X$.

\paragraph{1.2}
Soit $M\in{T_1(\overline{\mathbb{R}_o})}$ majorant de $X\subset{T_1(\overline{\mathbb{R}_o})}$, $M=m_0+m_1\cdot{o}$.
\paragraph{1er cas} $m_0>b_0$ alors $M<b_0+b_1\cdot{o}$.
\paragraph{2\`{e}me cas} $m_0=b_0$. $M=b_0+m_1\cdot{o}$ est un majorant de $X\cap{]b_0[}$ donc $b_1\leq{m_1}$ puisque $b_1$ est la borne sup\'{e}rieure de $ImF$.
\\Dans les deux cas, $b_0+b_1\cdot{o}\leq{M}$.

\paragraph{2} Que se passe-t-il en plus lorsque $b_I=-\infty$ ou $b_I=+\infty$ ?
\\Dans la d\'{e}monstration de $P(2)$, on montre que si $X$ est major\'{e} et $T_1(X)$ born\'{e} sup\'{e}rieurement par $b_0+\epsilon$ ou $b_0-\epsilon$, on a aussi $T_2(X)$ born\'{e} sup\'{e}rieurement par $b_0+\epsilon$ ou $b_0-\epsilon$.

\paragraph{2.a} Si $(\forall{x\in{X}})$ $T_1(x)\leq{b_0+\epsilon}$, alors $T_2(x)\leq{b_0+\epsilon}$ aussi et $b_0+\epsilon$ est majorant de $T_2(X)$ dans $T_2(\overline{\mathbb{R}_o})$. On montre que c'est le plus petit.
\\
Imaginons que $(\exists{m}\in{T_2}(\overline{\mathbb{R}_o}))$ $(\forall{x\in{X}})$ $T_2(x)\leq{m}<b_0+\epsilon$. Alors, $T_1(x)\leq{T_1(m)}$ et $T_1(m)\geq{b_0+\epsilon}$ puisque $b_0+\epsilon=SupT_1(X)$ dans $T_1(\overline{\mathbb{R}_o})$. On a $T_0(m)>b_0$ ou $T_1(m)=b_0+\epsilon$.
\\
Dans les deux cas $m\geq{b_0+\epsilon}$. Contradiction.

\paragraph{2.b} Si $(\forall{x\in{X}})$ $T_1(x)\leq{b_0-\epsilon}$, alors $T_0(x)<b_0$ ou $T_1(x)=b_0-\epsilon$. Dans les deux cas, $T_2(x)\leq{b_0-\epsilon}$. C'est majorant de $T_2(X)$ dans $T_2(\overline{\mathbb{R}_o})$, on montre que c'est le plus petit.
\\
Imaginons que $(\exists{m}\in{T_2(\overline{\mathbb{R}_o})})$ $(\forall{x\in{X}})$ $T_2(x)\leq{m}<b_0-\epsilon$. Donc $T_1(m)$ est un majorant de $T_1(X)$ dans $T_1(\overline{\mathbb{R}_o})$ et $T_1(m)\geq{b_0-\epsilon}$ puisque
$b_0-\epsilon=SupT_1(X)$ dans $T_1(\overline{\mathbb{R}_o})$.
\\
On en d\'{e}duit que $m\geq{b_0-\epsilon}$ car il n'existe aucun nombre de partie standard $b_0$ inf\'{e}rieur \`{a} $b_0-\epsilon$. Contradiction.

\paragraph{2.c}
D\`{e}s que la borne sup\'{e}rieure de $X$ est l'extr\'{e}mit\'{e} infinie d'une coupure infinit\'{e}simale, elle le reste pour tous les $T_I(X)$ suivants, dans la topologie de $T_I(\overline{\mathbb{R}_o})$. On montre que c'est encore le cas pour $X$ dans la topologie de $\overline{\mathbb{R}_o}$.
\\
On a $(\forall{i\geq{I}})$ $SupT_i(X)=b=b_0+b_1\cdot{o}+\dots{+b_{I-1}}\cdot{o^{I-1}}+\infty\cdot{o^{I}}$ (resp. $-\infty\cdot{o^{I}}$) dans $T_i(\overline{\mathbb{R}_o})$ et $(\forall{i<I})$ $SupT_i(X)=T_i(b)$ dans $T_i(\overline{\mathbb{R}_o})$. On montre que $Sup(X)=b$ dans la topologie de $\overline{\mathbb{R}_o}$.

\paragraph{2.c.1}
Si $(\exists{x}\in{X})$ $x>b$, $T_{I-1}(X)>T_{I-1}(b)$ (resp. $T_{I}(x)>b$).
\\Ceci est contradictoire avec $SupT_{I-1}(X)=T_{I-1}(b)$ dans $T_{I-1}(\overline{\mathbb{R}_o})$ (resp. $SupT_{I}(X)=b$ dans $T_{I}(\overline{\mathbb{R}_o})$).

\paragraph{2.c.2}
On montre que $b$ est le plus petit majorant de $X$.
\\Si $m\in{\overline{\mathbb{R}_o}}$ est un majorant de $X$ tel que $m<b$, alors $T_I(m)<b$ (resp. $T_{I-1}(m)<T_{I-1}(b)$) ce qui est contradictoire avec $SupT_{I}(X)=b$ dans $T_{I}(\overline{\mathbb{R}_o})$ (resp. $SupT_{I-1}(X)=T_{I-1}(b)$ dans $T_{I-1}(\overline{\mathbb{R}_o})$).

\paragraph{3}
On consid\`{e}re maintenant le cas o\`{u} tous les $b_i$ sont finis et l'on montre que $b=\sum\limits_{i\geq{0}}b_io^{i}=Sup(X)$ dans $\overline{\mathbb{R}_o}$.
\\
$X$ est major\'{e} et $(\forall{i}\in{\mathbb{N}})$ $SupT_i(X)=T_i(b)$ dans la topologie de $T_i(\overline{\mathbb{R}_o})$. Deux nouveaux raisonnements par l'absurde suffisent.

\paragraph{3.1}
Si $x>b$, il existe $i_0$ tel que $T_{i_0}(x)>T_{i_0}(b)$ et $T_i(x)=T_i(b)$ pour $0\leq{i}\leq{i_0}$ ce qui est contradictoire avec l'hypoth\`{e}se au rang $i_0$. On montre que $b$ est le plus petit majorant de $X$.

\paragraph{3.2}
Si $m\in\overline{\mathbb{R}_o}$ est un majorant de $X$ et $m<b$, on a aussi $i_0$ tel que $T_{i_0}(m)<T_{i_0}(b)$ et $T_i(m)=T_i(b)$ pour $0\leq{i}\leq{i_0}$ or $T_{i_0}(m)$ est un majorant de $T_{i_0}(X)$ et donc, puisque $SupT_{i_0}(X)=T_{i_0}(b)$ dans $T_{i_0}(\overline{\mathbb{R}_o})$ on a $T_{i_0}(b)\leq{T_{i_0}(m)}$. Contradiction. }
\qed

\begin{rem} $\overline{\mathbb{R}_o}$ est un dictionnaire de num\'{e}ros car il est totalement ordonn\'{e} et, comme un dictionnaire de lettres, cela n'aurait aucun sens d'ajouter des mots.
\par
Un dictionnaire fini n'est pas dense puisqu'il y a la notion de successeur. Un dictionnaire infini de tous les mots form\'{e}s d'un nombre fini de lettres serait dense mais pas continu.
\par
Si l'on consid\'{e}re tous les mots ayant un nombre fini ou infini d\'{e}nombrable de lettres, le "dictionnaire" devient continu car l'alphabet est discret mais ce n'est pas le cas de $\mathbb{R}$ (il n'y a pas de successeur dans $\mathbb{R}$).
\end{rem}

\begin{prop}
$(\Omega,\leq)$ n'est pas continu.
\end{prop}
Preuve : elle est la m\^{e}me que pour la \textbf{Proposition 5.18}. \qed
\\ \\
\begin{defn} On pose $\overline{\Omega}=\Omega\cup\mathbb{F}$ avec
$$\mathbb{F}=\{\sum\limits_{1\leq{k}\leq{N}}
b_k\Sigma^k+t+\sum\limits_{1\leq{p}\leq{P}}a_p\cdot{o^p}/N\in\mathbb{N}^*,
b_k,t\in\mathbb{R},$$
$$P\in\mathbb{N}^*
,a_P=\pm\infty,a_p\in\mathbb{R}, 1\leq{p}<P\}.$$
\end{defn}

\begin{thm}
$(\overline{\Omega},\leq)$ est continu.
\end{thm}
Preuve : $\overline{\Omega}$ est dense dans lui-m\^{e}me, on montre qu'il v\'{e}rifie le Th\'{e}or\`{e}me de Bolzano.
Soit $X\subset{\overline{\Omega}}$ non vide et minor\'{e}. L'ensemble de ses parties enti\`{e}res est not\'{e} $[X]$ et $[X]\subset{\aleph}$ est non vide et minor\'{e} par $M$.
\\
Il existe $L_0\in{\aleph}$ tel que $L_0$ soit un minorant de $[X]$ et $L_0+1$ ne le soit plus car si ce n'\'{e}tait pas le cas on aurait, pour tous les $L\in{\aleph}$, si $L$ est un minorant de $[X]$ alors $L+1$ l'est aussi.
Puisque $M$ est un minorant de $[X]$, on prouverait par induction transfinie dans $\aleph$ que tous les entiers $L\geq{M}$ sont des minorants de $[X]$ et $[X]$ serait alors vide. Contradiction.
\\ \\
On peut affirmer alors que $L_0-1$ est un minorant de $X\subset{\overline{\Omega}}$ et que $L_0+1$ ne l'est plus. On utilise la translation $T:\overline{\Omega}\rightarrow{\overline{\Omega}}$, $x\mapsto{x-L_0}$ pour prouver que
-1 est un minorant de $T(X)$ mais pas 1. On d\'{e}montre comme pour le \textbf{Th\'{e}or\`{e}me 5.24} que $T(X)$ a une borne inf\'{e}rieure $b\in{[-1,1]}$ et $X$ a une borne inf\'{e}rieure unique $L_0+b$.
\qed

\subsection{Une nouvelle caract\'{e}risation des structures \hbox{$(\Omega,+,\times,\leq)$} et  \hbox{$(\overline{\Omega},\leq)$}}
$\mathbb{R}(\Sigma)=\mathbb{R}(o)$ est le corps des fractions de l'anneau des polyn\^{o}mes de degr\'{e} fini $\mathbb{R}[\Sigma]$ ou $\mathbb{R}[o]$. C'est le plus petit sur-corps de $\mathbb{R}$ contenant $\Sigma$ et son inverse $o$ (cf. la \textbf{Remarque 1.5}).
\\
On le munit de l'ordre lexicographique tel que
$$0<o\ll{1}\ll{\Sigma}.$$

\begin{rem} Cette structure alg\'{e}brique et ordinale n'est pas \emph{continue} car la coupure form\'{e}e de
\begin{center}
$C_{g}=\{\frac{P(o)}{Q(o)}\in\mathbb{R}(o)^+/(\frac{P(o)}{Q(o)})^2\leq{1+o}\}\cup{\mathbb{R}(o)^-}$
\end{center}
 et
\begin{center}
$C_{d}=\{\frac{P(o)}{Q(o)}\in\mathbb{R}(o)^+/(\frac{P(o)}{Q(o)})^2\geq{1+o}\}$
\end{center}
d\'{e}finit une s\'{e}rie formelle $(1+o)^{\frac{1}{2}}=1+\frac{o}{2}-\frac{o^2}{8}+\frac{o^3}{16}-5\frac{o^4}{128}\dots$
qui n'est pas une fraction rationnelle. Elle n'est pas non plus \emph{compl\`{e}te} car les troncatures successives de cette "racine" de $1+o$ forment une suite de Cauchy qui n'est pas convergente dans $\mathbb{R}(o)$.
\end{rem}
$\Omega$ (resp. $\overline{\Omega}$) \'{e}tant un corps (resp. ensemble) totalement ordonn\'{e}, on va d\'{e}montrer que $\Omega$ est le compl\'{e}t\'{e} de $\mathbb{R}(\Sigma)=\mathbb{R}(o)$ et que $\mathbb{R}(\Sigma)=\mathbb{R}(o)$ est continu dans $\overline{\Omega}$. On pourra alors dire :

\begin{center}
\textbf{R\'{e}sultats principaux}
\\$(\Omega,+,\times,\leq)$ est le \emph{plus petit} sur-corps strict de $\mathbb{R}$
 \\totalement ordonn\'{e} et complet.
 \\$(\overline{\Omega},\leq)$ est le plus petit sur-ensemble strict de $\mathbb{R}$
  \\totalement ordonn\'{e} et continu.
\end{center}

\begin{lem} En tant qu'alg\`{e}bres, $\mathbb{R}(o)\subset{\Omega}$.
\end{lem}
Preuve : toute fraction rationnelle $\frac{P(o)}{Q(o)}$ est d\'{e}veloppable en \emph{une} s\'{e}rie formelle $S(o)\in\mathbb{R}_{o}=\mathbb{R}[[o]]$ ssi $Q(o)$ n'est pas infinit\'{e}simal ($\frac{P}{Q}$ n'admet pas 0 pour p\^{o}le [8, p.240]).
\\
Si ce n'est pas le cas, $Q(o)=o^{k}Q_{1}(o)$ avec $\frac{P(o)}{Q_{1}(o)}=S_{1}(o)$. Alors, $\frac{P(o)}{Q(o)}=\Sigma^{k}\times{S_1(o)}\in\Omega$. \qed

\begin{defn} On note "s\'{e}rie formelle", un \'{e}l\'{e}ment quelconque de $\Omega$.
\end{defn}

\begin{thm} $\Omega$ est l'ensemble des limites des suites de Cauchy d'\'{e}l\'{e}ments de $\mathbb{R}(o)$ muni de son ordre lexicographique. Autrement dit, $\Omega$ est le compl\'{e}t\'{e} de $\mathbb{R}(o)$.
\end{thm}
Preuve : elle est imm\'{e}diate. Soit une suite de Cauchy $(\frac{P_n(o)}{Q_n(o)})_{n\in\mathbb{N}}$ de $\mathbb{R}(o)\subset\Omega$. On sait que $\Omega$ est complet (cf. la  \textbf{Proposition 5.12}), donc la suite converge vers un \'{e}l\'{e}ment unique de $\Omega$.
\\
R\'{e}ciproquement, soit $S(o)$ une "s\'{e}rie formelle" de $o$ qui n'est pas le quotient exact de deux polyn\^{o}mes de $o$ par la division selon les puissances croissantes. Les troncatures successives de cette "s\'{e}rie" forment une suite de Cauchy de $\mathbb{R}(o)$ qui converge \'{e}videmment vers cet \'{e}l\'{e}ment de $\Omega$. \qed

\begin{lem} $\mathbb{R}(o)$ est dense dans $\Omega$ et $\Omega$ est dense dans $\overline{\Omega}$.
\end{lem}
Preuve : soient $S(o)$ et $T(o)$ deux "s\'{e}ries formelles" telles que $S(o)\leq{T(o)}$. Il existe au moins une fraction rationnelle strictement comprise entre $S(o)$ et $T(o)$, c'est la somme finie $F(o)=T_N(\frac{S_d(o)+S_g(o)}{2})$ avec $N=ord(S_d(o)-S_g(o)))$ qui est un nombre fini puisque $S_d(o)\neq{S_g(o)}$.
\\ \\
On d\'{e}montre que $\mathbb{R}_o$ est dense dans $\overline{\mathbb{R}_o}$.
\\
Soient $x=t+\sum\limits_{k\geq{1}}a_ko^k$ (si $a_N=\pm\infty$ alors $a_k=0$ pour $k>N$) et $x'=t'+\sum\limits_{k\geq{1}}b_ko^k$ (si $b_P=\pm\infty$ alors $b_k=0$ pour $k>P$).
\par
On suppose que $x<x'$. Si $x,x'\in\mathbb{R}_o$, $x<\frac{x+x'}{2}<x'$. Il est facile dans tous les autres cas ($N\neq{P}$, $N=P$ et $x$ et $y$ dans une m\^{e}me coupure ou non) de trouver un polyn\^{o}me de $o$ strictement compris entre $x$ et $x'$.
\\
On se ram\`{e}ne \`{a} ce cas (comme dans la d\'{e}monstration du \textbf{Th\'{e}or\`{e}me 5.28}) pour d\'{e}montrer que $\Omega$ est dense dans $\overline{\Omega}$. \qed

\begin{rem} On sait que $\mathbb{R}(o)$ est discontinu (cf. la \textbf{Remarque 5.29}). $\mathbb{R}(o)$ n'est pas continu dans $\Omega$ car la m\^{e}me coupure que pr\'{e}c\'{e}demment, avec $o$ au lieu de $1+o$ :
 \begin{center}
$C_{g}=\{\frac{P(o)}{Q(o)}\in\mathbb{R}(o)^+/(\frac{P(o)}{Q(o)})^2\leq{o}\}\cup{\mathbb{R}(o)^-}$
\end{center}
 et
\begin{center}
$C_{d}=\{\frac{P(o)}{Q(o)}\in\mathbb{R}(o)^+/(\frac{P(o)}{Q(o)})^2\geq{o}\}$
\end{center}

donne
\begin{center}
$C_g=\mathbb{R}(o)^-\bigcup{([0[\cap{\mathbb{R}(o)})}$ et $C_d=\mathbb{R}(o)^+\backslash{([0[\cap{\mathbb{R}(o)})}$
\end{center}
 et $+\infty\cdot{o}=\epsilon$ n'appartient pas \`{a} $\Omega$ mais \`{a} $\overline{\Omega}$.
\end{rem}

\begin{thm} On peut assimiler $\overline{\Omega}$ \`{a} l'ensemble des coupures de $\mathbb{R}(o)$. Autrement dit, $\mathbb{R}(o)$ est continu dans $\overline{\Omega}$.
\end{thm}
Preuve : soit $(C_g,C_d)$ une coupure quelconque de $\mathbb{R}(o)$. Soient
$$C'_d=\{S_d(o)\in\overline{\Omega}/(\forall{F_g(o)\in{C_g}})F_g(o)
\leq{S_d(o)}\},$$  $$C'_g=\{S_g(o)\in\overline{\Omega}/(\forall{F_d(o)\in{C_d}})S_g(o)
\leq{F_d(o)}\}.$$
On d\'{e}montre que $(C'_g,C'_d)$ est une coupure de $\overline{\Omega}$.
\\ \\
On a $C_g\subset{C'_g}$ et $C_d\subset{C'_d}$. Soit $S(o)\in\overline{\Omega}$.
\\Deux cas sont \`{a} consid\'{e}rer :
\\$(\exists{F_d(o)\in{C_d}})F_d(o)<S(o)$ alors $S(o)\in{C'_d}$ et
\\
$(\forall{F_d(o)\in{C_d}}) S(o)\leq{F_d(o)}$ et $S(o)\in{C'_g}$.
\\
$\overline{\Omega}\subset{C'_g\cup{C'_d}}$ et puisque $C'_g\cup{C'_d}\subset{\overline{\Omega}}$, on a la premi\`{e}re condition.
\\
Soient $S_d(o)$ un \'{e}l\'{e}ment quelconque de $C'_d$ et $S_g(o)$ un \'{e}l\'{e}ment quelconque de $C'_g$. Il faut montrer que $S_g(o)\leq{S_d(o)}$.
\par
Si ce n'\'{e}tait pas le cas, on aurait $S_d(o)<S_g(o)$ puisque $\overline{\Omega}$ est totalement ordonn\'{e}. Il existerait  alors par le Lemme pr\'{e}c\'{e}dent une fraction rationnelle $F(o)$ qui serait strictement comprise entre $S_d(o)$ et $S_g(o)$.
\\
Il y a une contradiction, soit parce que $F(o)\in{C_g}$ et $S_d(o)<F(o)$, soit parce que $F(o)\in{C_d}$ et $F(o)<S_g(o)$.
\\ \\
Puisque $\overline{\Omega}$ est continu :
$$(\exists{!S(o)\in{\Omega}})(\forall{S_g(o)\in{C'_g}})(\forall
{S_d(o)\in{C'_d}}) S_g(o)\leq{S(o)}\leq{S_d(o)}.$$
En particulier :
 $$(\forall{F_g(o)\in{C_g}\subset{C'_g}})(\forall{F_d(o)\in{C_d}
 \subset{C'_d}}) F_g(o)\leq{S(o)}\leq{F_d(o)}.$$
$S(o)$ est unique \`{a} partager la coupure $(C_g,C_d)$ car s'il existait $T(o)$ ayant la m\^{e}me propri\'{e}t\'{e}, le m\^{e}me Lemme prouverait qu'il existe une fraction rationnelle entre $S(o)$ et $T(o)$. \qed
\\ \\
$\mathbb{R}(o)$ est donc continu dans $\overline{\Omega}$. Comme $\mathbb{R}(o)$ est le plus petit sur-corps de $\mathbb{R}$ contenant l'\'{e}l\'{e}ment nouveau $o$, on peut "dire" que $\overline{\Omega}$ est le plus petit sur-ensemble de $\mathbb{R}$ qui soit continu et totalement ordonn\'{e}.

\section*{Conclusion}
Outre une d\'{e}monstration enti\`{e}rement nouvelle du Th\'{e}or\`{e}me Fondamental de l'Analyse, les r\'{e}sultats principaux de cette recherche portent sur les ensembles de nombres :
\\ \\
\textbf{1. Le mod\`{e}le non standard de l'Arithm\'{e}tique de Peano \hbox{$(\aleph^+,+,\times,\leq,S)$} est "le" prolongement \emph{intrins\`{e}que} de l'ensemble standard $\mathbb{N}$ car il ne d\'{e}pend d'aucun param\`{e}tre choisi arbitrairement.}
\\ \\
\textbf{2. $(\mathbb{R}_o,+,\times,\leq)$ est "le" plus petit sur-anneau int\`{e}gre de $\mathbb{R}$ totalement ordonn\'{e}, intris\`{e}que et \emph{complet} mais il n'est pas continu.
\\
$(\overline{\mathbb{R}_o},\leq)$ est le plus petit sur-ensemble de $\mathbb{R}$ intrins\`{e}que et \emph{continu}.}
\\ \\
\textbf{3. $(\Omega,+,\times,\leq)$ est le plus petit sur-corps de $\mathbb{R}$ totalement ordonn\'{e}, intrins\`{e}que et \emph{complet} mais il n'est pas continu.
\\
$(\overline{\Omega},\leq)$ est le prolongement par continuit\'{e} du plus petit sur-corps $\mathbb{R}(o)$ de $\mathbb{R}$, totalement ordonn\'{e} par l'ordre lexicographique. Il est intrins\`{e}que, i.e. ses propri\'{e}t\'{e}s sont toutes "n\'{e}cessaires"}
\\ \\
Autrement dit, les ensembles $\mathbb{R}_o$, $\overline{\mathbb{R}_o}$, $\Omega$, $\overline\Omega$ et $\aleph^+$ sont les plus "simples" extensions \emph{intrins\`{e}ques} des ensembles $\mathbb{R}$ d'une part, $\mathbb{N}$ d'autre part. \\ \\
Dans un article prochainement d\'{e}pos\'{e} dans Archiv D.S. [5] on donne les propri\'{e}t\'{e}s tr\`{e}s voisines des ensembles $\mathbb{R}_3=\frac{\mathbb{R}[X]}{(X^3)}$ et $\aleph_3$ mais ce sont des prolongements \emph{non intrins\`{e}ques} de $\mathbb{R}$ et $\mathbb{N}$ puisqu'ils d\'{e}pendent du choix d'un param\`{e}tre $I=3$.

\par Ces ensembles de nombres non standard permettent une \'{e}tude enti\`{e}rement originale du probl\`{e}me de la gravitation \`{a} deux corps.

\section*{Remerciements}
Le 7 Janvier 2011, j'ai fait un s\'{e}minaire au Laboratoire M.A.M. de l'U.B.S. (Vannes) o\`{u} j'ai pr\'{e}sent\'{e} une premi\`{e}re version des deux premi\`{e}res parties de cet article, mais pour l'ensemble de nombres r\'{e}els $\mathbb{R}_3=\mathbb{R}[X]/(X^3)$. Dans une troisi\`{e}me partie, je pr\'{e}sentais seulement "deux jolis r\'{e}sultats \`{a} v\'{e}rifier soigneusement" sur les ensembles $\aleph$ et $\Omega$ tels qu'ils ont \'{e}t\'{e} ici d\'{e}finis.
\par
Je remercie Fr\'{e}d\'{e}ric Math\'{e}us de s'y \^{e}tre int\'{e}ress\'{e} et de m'avoir encourag\'{e} \`{a} publier ces "jolis r\'{e}sultats" sous ArXiv.
\par
Tous mes remerciements \'{e}galement \`{a} Qiyu Jin (M.A.M.) de m'avoir aid\'{e} \`{a} ma\^{\i}triser le Langage LateX.

\end{document}